\documentclass[12pt]{article}
\usepackage{amssymb, fullpage}

\def\be{\begin{equation}}
\def\ee{\end{equation}}

\def\R{{\sf I\kern-.2em R}}
\def\N{{\sf I\kern-.2em N}}
\def\C{\kern.1em{\raise.47ex\hbox{$\scriptscriptstyle
$}}\kern-.40em{\sf C}}
\def\Z{{\sf Z\kern-.32em Z}}

\def\be{\begin{equation}}
\def\ee{\end{equation}}

\newtheorem{theorem}{\noindent Theorem}
\newtheorem{lemma}{\noindent Lemma}
\newtheorem{definition}{\noindent Definition}
\newtheorem{corollary}{\noindent Corollary}
\newtheorem{statement}{\noindent Proposition}

\def\Iso{\mathop{\rm Iso}}
\def\Mat{\mathop{\rm Mat}}
\def\conv{\mathop{\rm conv}}
\def\ext{\mathop{\rm ext}}
\def\Cl{\mathop{\rm Cl}}
\def\span{\mathop{\rm span}}

\title{Random metric spaces and universality}
\author{A.~M.~Vershik\thanks{%
St.~Petersburg Department of Steklov Institute of Mathematics,
27 Fontanka, 191011 St.~Petersburg, Russia.
Partially supported by the RFBR grant 02-01-00093 and the President
of Russian Federation grant for support of leading scientific schools 
2251.2003.1.}}
\date{}

\begin{document}
\maketitle
\abstract
{We define the notion of a random metric space and prove that
with probability one such a space is isometric to the Urysohn
universal metric space. The main technique is the study of
universal and random {\it distance matrices}; we relate the properties of
metric (in particular, universal) spaces to the properties
of distance matrices. We give examples of other categories
in which the randomness and universality coincide (graphs, etc.).}

\begin{center}
CONTENTS
\end{center}

1. Introduction. Relation between randomness and universality.

2. The cone of distance matrices and its geometry. Random distance matrices.

3. Universal distance matrices and the Urysohn space.
Urysohnness criterion.

4. The main theorem: universality of almost all 
distance matrices.

5. Matrix distributions as a complete invariant of metric triples
and a generalization of Kolmogorov's problem on extension
of measures.

Appendix. Examples of categories where universality and randomness 
coincide. The universal and random Rad\'o--Erd\"os--Renyi graph, etc.

\section{Introduction. Relation between randomness and universality}

Random objects in various concrete, relatively simple algebraic
and geometric categories (random
graphs, random polyhedra, etc.)
have been studied for a long time; however, as was already
mentioned in \cite{V4}, only recently such studies
started to assume systematic character; one may conjecture that in the near future 
this direction will become one of the most actual research fields.
This is motivated by numerous applications of
stochastic models, in particular in the field theory and statistical physics,
as well as by needs of mathematics itself, for example in the 
representation theory, number theory, theory of dynamical systems, and
the probability theory itself.

The classical probability theory and theory of random processes 
dealt and deal mainly with scalar or vector random variables, maybe
infinite-dimensional. But the foundations of probability theory,
which were formulated in the famous work by A.~N.~Kolmogorov on the basis
of general measure theory, of course, make it possible to develop the 
theory of random objects in arbitrary categories and formulate laws 
of large numbers and limit theorems that are quite far from the classical
ones. For example, we may speak of the limit distribution of
partitions of positive integers and Young diagrams, the
``law of large representations'' for random representations in 
series of classical groups (see papers on the asymptotic representation
theory, e.g., in \cite{Spr}), or the ``laws of large metric spaces''
(this will be done in this paper), etc. It is interesting that the technical
part of the corresponding considerations uses rather classical probabilistic
ideas and terms, but requires solving nontraditional problems.
It is equally important that these problems lead to new remarkable 
distributions, which did not arise earlier (for example, the Airy
distribution, the Poisson--Dirichlet distribution, the Lebesgue
measure in the space of subordinators, etc., which arose in the theory of
random partitions, Young diagrams, or spectra of random matrices).

Some time ago the author became interested in the following question:
what is a random metric space, and what can it look like if this
question makes any sense. This question arose in connection
with some problems of measure theory and ergodic theory (the study
of metrics related to classification of filtrations, i.e., decreasing
sequences of $\sigma$-algebras) from one hand, and 
with the study of the remarkable Urysohn universal metric space (see below), 
which had been  undeservedly forgotten and which manifestly pretended 
to play the central role in such considerations,
from the other hand. {\it  The main result of this
paper and the recent publications {\rm\cite{V2,V5}} is that the Urysohn space
is not only universal, but also generic in the topological (categorical)
sense, and, which is most important, it is ``random'' in the
following sense: ``sufficiently random''
finite metric spaces converge to the Urysohn space with probability one}.

Comparing with a simpler situation, namely with the theory of random graphs,
developed in the 60s, where the universal graph in the sense of Rad\'o
turns out, in a trivial way, to be also random in the Erd\"os--R\'enyi sense
(and, in fact, in much stronger sense --- see Section 5),
allows one to believe that the coincidence of universality
and randomness is not random and occurs in many categories.
A brief survey of categories where such a coincidence takes place 
is given in Appendix and will be considered in more detail elsewhere.

There are several fields related to the problems discussed in this paper;
we will mention them only briefly and in connection with concrete problems.
First, this is the {\it theory of random  matrices}, since we consider
random distance matrices. A special class of such matrices ---
matrix distributions --- are complete invariants of metric triples (or
Gromov triples, or $mm$-spaces), i.e., Polish spaces with Borel probability
measure, regarded up to measure-preserving isometries. A
matrix distribution is in turn the distribution of a random distance matrix,
which is invariant and ergodic under the action of the group
of simultaneous permutations of rows and columns and satisfies a certain
technical condition (the simplicity condition, see Section 5); these
relations and arising new problems will be mentioned below. In particular,
it is especially interesting to study the limit spectrum of the random
distance matrix, or, in other words, the random metric on the set of
positive integers. Second, the universality problems under consideration, 
as far as we deal with countable spaces, are related to logic, namely to the
    {\it theory of models } (see, e.g., \cite{P}), where 
    one considers the universality
(categoricity of a model) for different theories, and  
in particular for graph theory. For example, there is a striking 
relation between $0-1$ laws in the theory of first-order languages 
(see \cite{Gl, Bu, Hr, KC, Ca, Fa}) and similar 
Kolmogorov-type laws in the theory of random processes.
Another parallel with the classical probability theory 
is related to Kolmogorov's
theorem on extension of a cylinder measure to a countably additive
measure in the infinite product of spaces; we consider its nonlinear
analog: 
the problem of extending a finitely additive measure
defined on the set algebra generated by the family of balls in a metric
space centered at the points of a countable everywhere dense system;
it turns out that such an extension exists not for all spaces
(see Section 5).

Let us turn to a more detailed description of the main results and contents
of the paper. We begin with a description of the class (category)
of equipped metric spaces, i.e., complete separable (= Polish)
metric spaces with a distinguished countable well-ordered everywhere dense 
subset. We identify an equipped Polish space with the matrix formed
by the distances between the points of the distinguished subset;
this matrix is in a certain sense the set of structural constants
of the metric space. Such a viewpoint allows us to introduce a topology,
define measures, etc., on the space (category) of equipped metric spaces. 

In the second section we study the geometry of the {\it cone of
distance matrices} and the  convex set consisting 
of the vectors formed by the distances from the
points of a given finite metric space to a new point 
attached to this metric space.
We define 
the procedure of {\it random growth of finite metric spaces}
and give a precise statement of the main problem. The machinery of
distance matrices is convenient for studying some properties of metric
spaces; for example, in the recent paper \cite{Sm},
it is shown how one can find the \v{C}ech homologies of a metric space,
given the distance matrix of this space. Recall that in the same time the classification
of Polish (noncompact) spaces up to isometry is a ``wild'' problem (see
\cite{K}). The main problem is to define the notion of random metric space
and to study the relation between randomness and universality. Note that every metric
triple (= Polish space with a nondegenerate Borel probability measure)
naturally generates a random distance matrix. Namely, this is the random
matrix of the distances between the points of countable (dense) subsets of
a realization of the sequence of i.i.d.~points of the space distributed 
according to this measure (see Appendix and \cite{G,V1,V2}). 
We have the following important
question: what invariants of the metric space are contained in the
spectrum of these random symmetric matrices; only
preliminary experiments related to this question
are carried out (see \cite{BB}). Recall that the
distribution of these matrices (the matrix distribution) is a complete
invariant of the metric triple (\cite{G,V1,V3}); hence the old question
``whether one can hear the form of the drum'' (i.e., recover a Riemannian
metric from the spectrum of the Laplacian) takes a new form: whether one
can recover a metric triple up to volume-preserving isometry from
the spectrum of the random matrices of the distances between independently
chosen random points of the space. Although the answer to the latter question
in the general case is most probably negative (as in the case
of the former question),
the random spectrum undoubtedly contains information on many invariants.
This question needs a detailed study.

In the third section we recall the definition of the universal space
$\Iso(\cal U)$, introduced by P.~S.~Urysohn in 1924 in the papers \cite{U,U1},
which were published already after his tragic death.
We introduce the notion of a universal distance matrix and then prove
that such and only such matrices can be the distance matrices of 
countable everywhere dense subsets of the Urysohn space. This 
``urysohnness'' criterion was also given in 
\cite{V2} and \cite{V5}; its proof generalizes and simplifies
the proof of the main assertions of the original paper by P.~S.~Urysohn 
on the existence and uniqueness of the universal Polish space.
Our criterion implies immediately that universal matrices form
an everywhere dense set in the cone of distance matrices and thus
implies the ``genericity'' of the Urysohn space within
the chosen model,
i.e., in the sense of the topology of the cone $\cal R$. It follows in turn
that with respect to a generic probability measure on the cone of distance
matrices, almost all matrices are universal and thus with probability one
determine spaces isometric to the Urysohn space. However, this
assertion does not allow us to find at least one such measure explicitly;
concrete examples are given in the next section. We give a survey
of known properties of the space $\cal U$.

In the fourth section we prove that the procedure defined in the
second section gives explicit examples of measures on the cone
of distance matrices with respect to which almost all distance matrices
are universal. Random distance matrices governed by these distributions
are the best, as far as allowed by the triangle inequality, simulation
of random matrices with independent entries. This gives another
justification of the thesis that the random metric space is universal.

In the fifth section we reproduce main facts on the classification of metric
triples up to measure-preserving isometry, which are described in
detail in \cite{V2,V3}. We state a theorem on the
{\it matrix distribution as a complete invariant of a metric triple}
and on the reconstruction of the metric triple
from this invariant and give a sketch of the proof.
This theorem allows us to study invariants of a metric space with measure
using only the matrix distribution. 
The following new effect is of particular interest: in metric spaces (for example, in 
the Urysohn space) there exist sequences that are uniformly distributed
with respect to a certain finitely additive measure, while there exists
no countably additive measure with respect to which they would be
uniformly distributed. Such examples are given by the same
construction from Section 3.

The thesis on coincidence of universality and randomness has
a much wider domain of applicability; in Appendix we recall the
combinatorial analog of our problem --- the relation between 
random and universal
graphs, as well as give other examples. Note that the papers \cite{Ra,ER,Ca} on
the universal graph appeared almost 40 years  later than Urysohn's
paper, but it was apparently unknown to their authors.
The example of graphs is a quite special case of our
scheme; namely, the universal graph is the universal metric
space in the class of metrics assuming only the values 
0, 1, 2 (hence in this case the triangle inequality is automatically satisfied).
But a deeper analogy between both theories is not only in the same
relation between universality and randomness, but in an analogy
between the group of isometries of the Urysohn space
and the group of automorphisms of the universal graph; these poorly studied
groups are of great interest as universal, in one or another sense, 
infinite-dimensional groups. The group of isometries of the {\it Urysohn
rational universal space}, which we will mention only briefly,
takes an intermediate position between these two examples.
The groups of symmetries of universal spaces in various categories
are of special interest from the purely algebraic point of view
as well as for the whole theory we develop.

The author is grateful to Erwin Schr\"odinger Institut (Vienna)
and Max-Planck-Institute (Bonn) for offering him the possibility to work. 
He is also grateful to
E.~Gordon, P.~Cameron, J.~Baldwin, Y.~Benjamini, B.~Granovsky, A.~Bovykin,
and V.~Uspensky for useful references obtained from them at various times. 
A part of the material of the paper
is made up from the lectures given by the author in 
Vienna University (Austria, October 2002), Santiago University (Chile,
December 2002), Les Houches Winter Physics School (France, March 2003),
and Haifa Technion (Israel, June 2003). The Ph.D.~student U.~Haboeck
helped me to process the texts of Vienna lectures.
The author also thanks the group of physicists 
(O.~Bohigos et al.) from Universit\'e
Orsay (Paris-IX), who became interested in these subjects after the
author's lecture in Oberwolfach and made very useful numerical 
experiments with spectra of distance matrices for simplest
manifolds.

\section{The cone $\cal R$ of distance matrices and its geometry.
Random distance matrices}

\subsection{Definitions}
Consider the following set of infinite real matrices:
 $${\cal R} =\{\{r_{i,j}\}_{i,j=1}^{\infty}:  r_{i,i}=0,\,r_{i,j}\geq 0,\,
 r_{i,j}=r_{j,i},\,  r_{i,k}+r_{k,j}\geq r_{i,j}{\rm \,\,for\,\,}
 i,j,k=1,2, \dots\}.$$
Elements of the set $\cal R$ will be called {\it distance matrices};
note that zeros outside the main diagonal are allowed. Every distance
matrix determines a semimetric on the set of positive integers
${\bf N}$. A matrix from $\cal R$ is called a {\it proper distance matrix}
and determines a metric on ${\bf N}$ if it has no zeros outside the
main diagonal. The set of all distance matrices is a convex cone in
the vector space of infinite real matrices
$\Mat_{\bf N}({\bf R})={\bf R^N}^2$; this cone is weakly closed 
in the ordinary weak topology in the space of infinite real matrices.
We will call it the {\it cone of distance matrices}.
The subset of proper distance matrices is an open weakly dense
subcone in ${\cal R}$.  If $r$ is a proper distance matrix, then
the completion of the metric space $({\bf N},r)$ is a complete
separable metric (= Polish) space ($X_r,\rho_r)$ equipped with
a countable everywhere dense well-ordered subset
$\{x_i\}_{i=1}^{\infty}$, which is the image of the set of 
positive integers in the completion. An arbitrary distance matrix
(with possible zeros outside the main diagonal) determines a semimetric
on the set of positive integers; in this case, before constructing
the completion, we should first take the quotient space with respect to 
the partition of $\bf N$ into classes of points with zero distance. For example,
the zero matrix leads to the one-point metric space.
Thus finite metric spaces are also included in our considerations. 

Now assume that we are given a Polish space
$(X,\rho)$ equipped with a countable everywhere dense well-ordered subset
$\{x_i\}_{i=1}^\infty$. The proper distance matrix
$r=\{r_{i,j}\} \in \cal R$, where $r_{i,j}=\rho(x_i,x_j)$,
$i,j =1,2,\dots$, determines a metric on the set of positive integers.
Thus the correspondence between proper distance matrices and equipped Polish
spaces is a bijection. Every invariant property of the metric space
(topological, homological, etc.) can be expressed in terms of the 
distance matrix of any countable everywhere dense set.

We may regard the cone $\cal R$ as a fiber bundle whose base space
is the class of all individual Polish spaces (note that,
by the universality of the Urysohn space (see below), we may
assume that the base space is the set of all closed subsets
of the Urysohn space) and the fiber over a given space is the set of all
countable well-ordered everywhere dense subsets of this space.
Thus the cone $\cal R$ is the {\it universum of equipped Polish spaces},
and we will study the properties of these spaces and the whole collection of them
with the help of this cone.

Similarly to the definition of the cone
$\cal R$, we can define the finite-dimensional cones 
${\cal R}_n$ of distance matrices of order $n$. The cone ${\cal
R}_n$ is a polyhedral cone lying inside the positive octant of the space
of matrices $ \Mat_n({\bf R})\equiv {\bf R}^{n^2}$.
Let us define a subset ${\bf M}^s_n({\bf R})\equiv {\bf M}^s_n$
of the space of symmetric matrices as the set of {\it symmetric matrices
with zeros on the main diagonal}. The cone ${\cal R}_n$
is contained in this set of matrices:
${\cal R}_n \subset {\bf M}^s_n$, and ${\bf M}^s_n$ is obviously
the linear span of the cone: $\span({\cal R}_n)={\bf
M}^s_n$, since the interior of
${\cal R}_n$ is not empty. Clearly,
$\span({\cal R})\subset {\bf M}^s_{\bf N}$, where ${\bf M}^s_{\bf
N}$ is the space of all real symmetric matrices with zeros on the
main diagonal; the geometry of the cone
$\cal R$ is quite nontrivial. Every distance matrix
$r \in {\cal R}_n$ determines a (semi)metric on the space
$X_r$ consisting of $n$ points.

Consider the projections
$$p_{m,n}:{\bf M}^s_m \longrightarrow {\bf M}^s_n,\quad m >n, $$
that associate with a matrix $r$ of order $m$ its NW-corner
(= north-west corner) of order $n$. The cones ${\cal R}_n$ agree with
these projections, i.e., $p_{m,n}({\cal R}_m)={\cal R}_n$.
The projections $p_{n,m}$ extend naturally to 
a projection $p_n: M^s_{\bf N}  \longrightarrow M^s_n({\bf R})$ on
the space of infinite symmetric matrices with zero diagonal,
and the resulting projection 
$p_n$ also preserves the cones:
$p_n({\cal R})={\cal R}_n$. It is clear that the cone ${\cal
R}$, regarded as a topological space (in the weak topology), is
the inverse limit of the sequence
$({\cal R}_n,\{p_n\})$.

An obvious but important property of the cones
${\cal R}_n$ is their invariance under the action of the
symmetric group $S_n$ whose elements permute simultaneously the rows
and columns of matrices.

Let us consider the geometric structure of the cones
${\cal R}_n$ and ${\cal R}$. In the first dimensions we have
${\cal R}_1=\{0\}$ and ${\cal
R}_2={\bf R}$. The description of extremal rays (in the sense of
convex geometry) of the convex polyhedral cone
${\cal R}_n$, $n=3,4,\dots, \infty$, is a well-known difficult
combinatorial-geometric problem (see \cite{DL, Av} and references
given there). Every extremal cone in
${\cal R}_n$, $n \leq 4$, is of the form $\{\lambda \cdot l: \lambda \geq
0\}$,  where $l$ is a symmetric matrix with entries
$0$ and $1$, i.e., the distance matrix associated with a semimetric
such that the corresponding quotient space
consists of two distinct points.
But if $n \geq 5$, there are extremal rays consisting
of proper metrics. A complete description of the set of extremal rays
for the general dimension is not known. The most interesting
question is the study of the
asymptotic properties of the cones
${\cal R}_n$, and especially the description of the extremal rays of the cone
${\cal R}$.         
It may happen that the set of these rays is a dense
$G_{\delta}$-set in ${\cal R}$ and some of them consist of universal 
distance matrices. This agrees with the estimates of the growth of
the number of extremal rays of the cone
${\cal R}_n$ given in  \cite{Av}. It seems that the 
algebro-geometric structure and
stratification of the cones 
${\cal R}_n$ regarded as semialgebraic sets was not studied.  
In order to study the topological and convex structure of the cones
${\cal R}_n$, we will use their 
{\it inductive description}, which is given below.

\subsection{Admissible vectors and the structure of the cone
${\cal R}_n$}

Let $r=\{r_{i,j}\}_1^n$ be a distance matrix of order
$n$ ($r\in {\cal R}_n$); choose a vector $a \equiv \{a_i\}_{i=1}^n \in
{\bf R}^n$ with the following property:
if we attach it to the matrix $r$ as the last row
and column, putting zero at the diagonal, we again obtain a distance matrix of order
$n+1$. We say that  such vectors are
{\it admissible} for the given fixed distance matrix $r$ 
and denote the set of admissible vectors for $r$
by $A(r)$. Let $(r^a)$ stand for the distance matrix of order
$n+1$ obtained from the matrix $r$ by attaching the vector $a \in
A(r)$ as the last row and column; obviously,
$p_n(r^a)=r$. The matrix  $r^a$ is of the form
 $$r^a = \left( \begin{array}{ccccc}
 0    &r_{1,2}& \ldots&  r_{1,n}&a_1\\
 r_{1,2} &0& \ldots &r_{2,n} &a_2 \\
 \vdots &  \vdots&
\ddots & \vdots &\vdots\\
r_{1,n}&r_{2,n}& \ldots &0  &a_n\\
  a_1&  a_2  & \ldots & a_n& 0
\end{array}\right).$$

The (semi)metric space $X_{r^a}$ corresponding to the matrix  $r^a$
is an extension of the space $X_r$: we add one new point
$x_{n+1}$, and the numbers $a_i$, $i=1, \dots, n$, are the distances from
$x_{n+1}$ to the old points $x_i$. The admissibility of the vector $a$
is equivalent to the following system of inequalities: an admissible vector
$a=\{a_i\}_{i=1}^n$ for a fixed matrix
$\{r_{i,j}\}_{i,j=1}^n$ must satisfy the series of triangle inequalities
for all $i,j,k=1,2,\dots, n$:
\begin{equation}
|a_i-a_j|\leq r_{i,j}\leq a_i+a_j.
\end{equation}

Thus for a distance matrix $r$ of order $n$, the set of
admissible vectors is $A(r)=\{\{a_i\}_{i=1}^n:|a_i-a_j|\leq
r_{i,j}\leq a_i+a_j,\; i,j=1, \dots, n\}$. It is worth mentioning that
a vector $a=\{a_i\}$ may be regarded as  a {\it Lipschitz function} $f_a(.)$,
$f_a(i)=a_i$,
on the space $X_r=\{1,2, \dots, n\}$ with metric $r$,
with Lipschitz constant  $1$. From this point of view it
is useful to consider extensions of metric spaces
(see \cite{Ka,G}).

Geometrically, the set of admissible vectors
$A(r)$ can be identified with the intersection of the cone
${\cal R}_{n+1}$ and the affine subspace that consists of matrices
of order $n+1$ whose 
NW-corner (north-west corner) 
coincides with the matrix $r$.
Since the inequalities are linear, it is clear that the set
$A(r)$ is an unbounded polyhedral convex set in
${\bf R}^n$. If $r_{i,j} \equiv 0$, $i,j=1, \dots, n$, $n\geq 1$, then $A(r)$
is the diagonal: $A(0)=\Delta_n \equiv \{(\lambda, \dots, \lambda):
\lambda \geq 0\}\subset {\bf R}^n_+$. Let us describe the structure of the 
sets $A(r)$ in more detail.

\begin{lemma}
{For every proper distance matrix $r$ of order $n$, the set of admissible 
vectors $A(r)$ is a closed polyhedral set lying in the octant
${\bf R}_+^n$, namely the Minkowski sum
$$A(r)= M_r +\Delta_n,$$
where $\Delta_n $ is the half-line of constant vectors in the space
${\bf R}_+^n$ and $M_r$ is a convex compact polytope
of dimension $n$, which is the convex hull of the set of extreme points
of the whole polyhedron $A(r)$: $M_r= \conv(\ext A(r))$.}
\end{lemma}

\begin{proof}
{The set $A(r) \subset {\bf R}^n$ is the intersection of finitely many
closed half-spaces and obviously does not contain lines. By general
theorems of convex geometry, the set $A(r)$ is the sum of a
convex polyhedron (polytope), which is the convex hull of the set of
extreme points of the set $A(r)$, and a closed polyhedral cone,
with vertex at the origin, that does not contain lines. It remains to observe that
this cone degenerates into the diagonal of the positive octant, i.e.,
into the half-line of vectors in ${\bf R}^n$ with equal coordinates;
indeed, if the cone contained another half-line, the triangle
inequality would break: the difference 
$a_i - a_j$ of the coordinates 
could be arbitrarily large for at least one pair
$(i,j)$. The dimension of the set
$A(r)$ equals $n$ for a proper distance matrix; in the general case
it depends on the matrix $r$ and can be less than $n$; 
the dimension of the polytope 
$M_r$ equals $\dim A(r)$ or $\dim
A(r)-1$.}
\end{proof}

The following lemma asserts that the correspondence
$r \rightarrow A(r)$ is covariant under the action of the
symmetric group in ${\bf R}^n$. Its proof is obvious.

\begin{lemma}
{For every matrix $r\in {\cal R}_n$, we have
$A(grg^{-1})=g(A(r))$, where $g \in S_n$ is an element of the symmetric
group $S_n$, which acts naturally on the space of matrices
$M_{\bf N}(\bf R)$ and on the space of convex subsets of the vector space
${\bf R}^n$}.
\end{lemma}

The convex structure of the polyhedral sets 
$M_r$ and $A(r)$ is very interesting; apparently it
has not been studied earlier. For dimensions greater than three, the 
{\it combinatorial type} of the polytope  $M_r$ 
essentially depends on $r$. But in the three-dimensional case, this combinatorial
type is the same for all proper distance matrices, and hence
the type of the set $A(r)$ is also the same for all proper distance matrices.
Let us consider this answer in more detail.

\smallskip\noindent
{\bf Example.} For $n=3$, the set $A(r)$ and the extreme points
of the polyhedron $M_r$ can be described  as follows. Let $r$ be a proper distance matrix
 $$r = \left( \begin{array}{ccc}
  0 & r_{1,2} & r_{1,3}\\
  r_{1,2}& 0 & r_{2,3} \\
  r_{1,3} &r_{2,3}&0 \\
\end{array}\right).$$
Set $r_{1,2}=\alpha$,  $r_{1,3}= \beta$,  $r_{2,3}=\gamma$;
then
 $$r^a = \left( \begin{array}{cccc}
 0    &\alpha  &  \beta & a_1\\
 \alpha &0& \gamma & a_2 \\
 \beta & \gamma &0  & a_3\\
  a_1&  a_2  &  a_3 &  0
\end{array}\right).$$
Let  $\delta=\frac{1}{2}(\alpha+\beta+\gamma)$. There are seven extreme points
$a=(a_1,a_2,a_3)$ of the polytope $A(r)$ : the first vertex is the closest
vertex to the origin: $(\delta-\gamma,
\delta-\beta, \delta-\alpha)$, three other ones are nondegenerate extreme
points (metrics): $(\delta,\delta-\alpha,
\delta-\gamma)$, $(\delta-\beta, \delta, \delta-\alpha)$,
$(\delta-\gamma, \delta-\beta, \delta)$, the remaining three vertices are 
degenerate: $ (0, \alpha, \beta)$, $(\alpha, 0, \gamma)$, $(\beta,
\gamma, 0)$.

If $\alpha=\beta=\gamma=1$, then these seven points are as follows:
$$(1/2,1/2,1/2), (3/2,1/2,1/2), (1/2,3/2,1/2), (1/2,1/2,3/2),
(0,1,1), (1,0,1), (1,1,0).$$

Note that all nondegenerate extreme points determine metrics 
on the space consisting of four points that cannot be 
isometrically embedded into an Euclidean space.

\subsection{Projections and isomorphisms}

Let $r$ be a distance matrix of order $N$, and let $p_n(r)$ 
be its NW- (north-west) corner of order $n<N$. Then we can define
a projection $\chi^r_n$ of the set $A(r)$ onto $A(p_n(r))$: $\chi^r_n:(b_1,
\dots, b_n, b_{n+1},\dots, b_N) \mapsto (b_1, \dots, b_n)$. (We omit
the index $N$ in the notation $\chi^r_n$.) The following simple lemma
plays a very important role in our further construction.

\begin{lemma} {\rm(}Amalgamation lemma{\rm)}.
{Let  $r\in {\cal R}_n$ be a distance matrix of order $n$. For any two
vectors $a=(a_1, \dots, a_n) \in A(r)$ and  $b=(b_1,
\dots, b_n) \in A(r)$, there exists a real nonnegative number
$h \in {\bf R}$ such that ${\bar b}=(b_1,\dots, b_n, h) \in
A(r^a)$ {\rm(}and also ${\bar a}=(a_1, \dots, a_n, h) \in A(r^b)${\rm)}.}
\end{lemma}

\begin{corollary}
{For every matrix $r \in {\cal R}_n$ and every vector $a \in A(r)$, the mapping
$\chi^r_{n+1,n}: (b_1, \dots, b_n, b_{n+1}) \mapsto (b_1, \dots,
b_n)$ of the set $A(r^a)$ into $A(r)$ is an epimorphism of the set
$A(r^a)$ onto $A(r)$
{\rm(}by definition, $p_{n+1,n}(r^a)=r${\rm)}}.
\end{corollary}
{\sloppy}

\begin{proof}
{The lemma follows from a simple geometric observation.
Consider two finite metric spaces:
$X=\{x_1, \dots,
x_{n-1}, x_n\}$ with metric $\rho_1$ and $Y=\{y_1, \dots, y_{n-1},y_n
\}$ with metric $\rho_2$. Assume that the subspaces consisting of 
the first $n-1$ points $\{x_1, \dots, x_{n-1}\}$ and
$\{y_1, \dots, y_{n-1}\}$ are isometric, that is,
$\rho_1(x_i,x_j)=\rho_2(y_i,y_j)$, $i,j=1, \dots, n-1$. Then there exists 
a third space $Z=\{z_1,\dots,
z_{n-1},z_n,z_{n+1}\}$ with metric $\rho$ and isometries $I_1$ and $I_2$
of the spaces $X$ and $Y$, respectively, into the space $Z$ such that
$$I_1(x_i)=z_i,\; I_2(y_i)=z_i,\; i=1, \dots, n-1,\; I_1(x_n)=z_n,\;
I_2(y_n)=z_{n+1}.
$$
In order to prove the existence of such a space 
$Z$, we must show that there exists a nonnegative number
$h$ that can be the distance 
$\rho(z_n,z_{n+1})=h$ between $z_n$ and $z_{n+1}$ (the images of $x_n$ and
$y_n$, respectively, in $Z$), i.e., such a distance $h$
must satisfy all triangle inequalities
in the space $Z$. The existence of such
a number $h$ follows from the inequalities
$$\rho_1(x_i,x_n)-\rho_2(y_i,y_n) \leq \rho_1(x_i,x_j)+
\rho_1(x_j,x_n)-\rho_2(y_i,y_n)$$
$$=\rho_1(x_j,x_n)+\rho_2(y_i,y_j)-\rho_2(y_i,y_n)\leq \rho_1(x_j,x_n)+
\rho_2(y_j,y_n),$$ 
which hold for all $i,j=1, \dots, n-1$.
Therefore, $$\max_i |\rho_1(x_i,x_n)-\rho_2(y_i,y_n)| \equiv M
\leq  m \equiv  \min_j (\rho_1(x_j,x_n)+\rho_2(y_j,y_n)).$$ 
Thus we can take $h$ equal to an arbitrary number from the
nonempty closed interval $[M,m]$, and we set $
\rho(z_n,z_{n+1}) \equiv h$; it follows from the definition of $h$ 
that all triangle inequalities are satisfied. Now assume that $r$
is a distance matrix of order $n-1$ and $a \in A(r)$ 
is an admissible vector; let $\{x_1, \dots, x_{n-1}, x_n\}$
be a space such that the distances between the first
$n-1$ points are given by the matrix
$r$ and the metric on the whole space is given by the extended matrix
$r^a$. We can choose another arbitrary admissible vector
$b \in A(r)$ and obtain the distance matrix
$r^b$ of a space $\{y_1, \dots, y_{n-1},y_n\}$, where the space
consisting of the first
$n-1$ points is isometric to the space
$\{x_1, \dots, x_{n-1}\}$. As we have proved above, there exists a space
$Z$ whose distance matrix 
$\bar r$ of order $n+1$ enjoys the desired property.}
\end{proof}

``Amalgamation lemmas,'' i.e., lemmas on extension of homomorphisms
of two spaces $X_1$ and $X_2$ to a space $X_3$ with an identified subspace
$X_4$ to a homomorphism of their union to $X_3$,
are an element of the construction of the inductive limit that
is known in the theory of models as the Fraisse limit. We apply the term
``amalgamation'' exactly in the same context in which
it is used  in the theory of models, when one studies the possibility
of embedding two intersecting but noncoinciding models
into a third one. In fact, this lemma predetermines  the existence of
the universal object in the category, in our case ---
the universal metric space. In the next section we will use the inequality
we have proved
in the inductive construction of measures.

Now we are able to state a general assertion on the projection
$\chi^r$.

\begin{lemma}
{For arbitrary positive integers $N$ and $n<N$ and every matrix $r
\in {\cal R}_N$, the mapping $\chi^r_n$ is an epimorphism of the set
$A(r)$ onto $A(p_n(r))$. In other words, for every
$a=(a_1, \dots, a_n) \in A(p_n(r))$, there exists a vector 
$(b_{n+1},\dots, b_N)$ such that $b=(a_1,\dots, a_n,  b_{n+1}, \dots, b_N)
\in A(r)$.}
\end{lemma}

\begin{proof}
{The proof of the previous lemma shows how we should determine the first 
number $b_{n+1}$. But $\chi^r_n$ is a mapping of the
set $A(r)$, $r\in {\cal R}_N$, to the similar set
$A(p_n(r))$ and is the product of the projections 
$\chi^r_n,\cdots, \chi^r_{N-1}$, each of them being an epimorphism; hence
the product is also epimorphic.}
\end{proof}

In what follows, it is convenient to represent an infinite matrix
 $r \equiv \{r_{i,j}\} \in {\cal R}$ as a sequence of admissible vectors
 of increasing lengths
\begin{equation}
r(1)=\{r_{1,2}\},\; r(2)=\{r_{1,3},r_{2,3}\}, \dots,
r(k)=\{r_{1,k+1},r_{2,k+1}, \dots, r_{k,k+1}\},\dots, \;k=1,2, \dots,
\end{equation}
satisfying the conditions  $r(k) \in A(p_k(r))$ (recall that
$p_k(r)$ is the NW-projection of a matrix
$r$ to the space ${\bf M}^s_k$, i.e., the NW-corner of $r$), 
that is, such that every vector
$r(k)$ is admissible for the
 {\it previous distance matrix $p_k(r)$}.

We may consider the following sequence of cones and mappings:
\begin{equation}
0={\cal R}_1\stackrel{p_2}{\longleftarrow} {\cal
R}_2={\bf R}_+\stackrel{p_3}{\longleftarrow}{\cal R}_3
{\longleftarrow} \dots {\longleftarrow}{\cal
R}_{n-1}\stackrel{p_n}{\longleftarrow}{\cal R}_n{\longleftarrow}
\dots;
\end{equation}
here $p_n$ is the restriction of the projection
defined above to the cone
${\cal R}_n$. The preimage of a point $r \in {\cal R}_{n-1}$ (the fiber over
$r$) is the set $A(r)$ defined in Lemma 1. Note that the sequence (3)
is not a sequence of fiber bundles in the ordinary sense: the preimages 
may be even nonhomeomorphic for different points $r$ (even
their dimensions may be different). But this sequence determines the cone
$\cal R$ as the inverse limit of the cones
${\cal R}_n$. We use the sequence 
(3) to define measures on the cone ${\cal R}$ in the spirit
of the theory of Markov processes.

\subsection{Definition of the ``most'' random metric space}

Using the geometry of the cone of distance matrices, we will define
a class of Borel probability measures on this cone that may be regarded
as the most ``dispersed'' distributions, as far as the 
cone of distance matrices allows.
A random distance matrix governed by this distribution corresponds
to the intuitive notion of a random metric space in the following
sense, which will be specified below: the set of supports of the 
distributions of finite fragments of this random matrix coincides with the
set of all distance matrices of finite metric spaces without any
restrictions. Since the cone ${\cal R}$ is infinite-dimensional, 
it has no distinguished measure (similar to the Lebesgue measure
on a finite-dimensional manifold); hence we should speak 
of the class of most
dispersed measures, rather than of 
a single measure.
This class, as we will see, includes measures
with very different behavior. On the other hand, in view of
the triangle inequality, we cannot (except for
the example described below)
assume that all pairwise distances are independent. We 
construct measures with the help of
a {\it Markovian procedure}, which imitates
independence and is useful in other cases. Briefly speaking,
we construct a random matrix as a random two-dimensional field,
Markovian in some sense.

The construction is inductive: we begin with the 
one-point space and then describe how to add a new
$(n+1)$th point to $n$ points already constructed, choosing the distances
from the new point to the old ones in the most random way. 
The ``maximum randomness'' means that, given a distance matrix
$r_n \in {\cal R}_n$ already chosen, we choose an admissible vector uniformly
distributed on $A(r_n)$. This method corresponds to the literal
understanding of what is a random finite metric space.

After infinitely many steps we obtain an infinite random distance matrix,
i.e., a probability measure on the cone
${\cal R}$, or, in other words, a random metric on the set of positive
integers. Then the following question arises: {\it what metric spaces
are the completions of the set of positive integers with respect
to these metrics}? We will answer this question in Section 3.
Now we mention only that this is {\it one and the same, up to isometry,
metric space} for almost all (in the sense of the measure
constructed) distance matrices, and almost all matrices are
the distance matrices of various everywhere dense countable
sequences of points in this space; this space
is precisely the universal Urysohn space.

\smallskip

Let us now describe the precise construction of probability 
measures on the cone
$\cal R$. The parameter of this construction is an arbitrary 
continuous probability 
measure $\gamma$ of full support
on the half-line ${\bf R}^1_+$,
for example the Gaussian measure on the half-line. Consider the one-point
metric space, and assume that its single point is the positive
integer $1$. Every time when we add a new point
$n$, $n= 2,3, \dots $, we let its distance
$r_{1,n}$ to the point
$1$ be equal to the random variable $\xi_n$,  $n = 2,3, \dots$, 
and assume that the random variables 
$\xi_n$ for different $n$ are i.i.d.~with distribution
$\gamma$. Thus the distribution of the first row (and the first column)
of the random matrix is defined; in particular, the distribution
of the entry $r_{1,2}$ is defined.

Now we define a measure $\nu_{\gamma}$ on the cone $\cal
R$ by induction. We use the following method of describing a measure:
fix $n>2$ and $k<n$, and assume that we are given a probability distribution
on the fragments of finite rectangular distance matrices, namely
on the set of collections $\{r_{i,j}\}$, where the indices range over
the sets $i,j=1, \dots, n-1$ and $(i,n)$, $i=1, \dots, k $;
in other words, we are given the joint distribution of the pairwise distances
between all pairs of the first
$n-1$ points and the distances from the $n$th point to the first $k$ points,
$1 < k < n$. Now let us fix the values of these distances and
determine the conditional distribution of the distance between the points $n$
and $k+1$.

Obviously, the unique inequality satisfied by 
the distance $r_{k+1,n}$ is
    $$  \max_{i=1,\dots, k} |r_{i,k+1}-r_{i,n}| \leq   r_{k+1,n}
    \leq \min_{i=1,\dots, k} |r_{i,k+1}+r_{i,n}|.$$
The amalgamation lemma from the previous section
implies that the left-hand side of the inequality does not exceed
the right-hand side, so that the range of possible values for
$r_{k+1,n}$ is not empty. {\it As the conditional measure for
$r_{k+1,n}$, take the Lebesgue measure on this interval normalized to one}.
Of course, the interval depends substantially on the matrix fragment we
consider; hence our measure is not (and, as mentioned above, cannot be for 
any change of parameters)
a product measure: the pairwise distances are not independent.
However, one should regard the above definition of conditional measures 
as the best admissible step-wise replacement of 
the independence of 
$r_{k+1,n}$ and the previous distances. When $k$ assumes the value 
$n-1$ (recall that $r_{n,n}=0$), the distribution of the matrix
$r_{i,j}$, $i,j=1,
\dots, n$, $n>2$, is defined, i.e., we obtain a random
$n$-point metric space. Since the distribution of the entry
$r_{1,n+1}$ is already determined, we define, according to the same rules,
the conditional measure of the entry
$r_{2,n+1}$, i.e., the distance between the points $2$ and 
$n+1$. The description of the algorithm is finished.

Our method of constructing a measure, which can be called Markovian,
is based on the following general fact from measure theory: if 
we are given a measurable partition (above it was the one-dimensional
fiber bundle over the described fragment of the matrix with fiber
$r_{k+1,n}$) in a standard measure space (Lebesgue--Rokhlin space),
then a measure on the whole space is uniquely determined  by its
projection to the base space (i.e., a measure on the quotient space)
and the system of conditional measures on fibers (for almost
all fibers). The fiber of our fiber bundle is the interval of
possible values for the distance between the $(k+1)$th and $n$th points,
given fixed (already determined) distances between the pairs of points
with numbers less than $n$ and between the $n$th and $i$th points for
$i<k+1$.

Thus the distribution of the entries of the first row and the above
construction of conditional measures uniquely determine the measure
$\nu_{\gamma}$.

\begin{lemma}(on stationarity).

{1. The support of the measure $\nu_{\gamma}$ is the whole cone $\cal R$.

2. Denote by $\nu_{\gamma}^n$ the projection of the measure $\nu_{\gamma}$
to the set of submatrices 
$\{\{a_{i,j}\}:\;i=1, \dots, n;\, j=n+1,n+2, \dots \}$, i.e., the joint
distribution of these matrix entries. The measure
$\nu_{\gamma}^n$ is invariant under the one-sided shift
$\{a_{i,j}\} \to \{a_{i,j+1}\}$ and ergodic.}
\end{lemma}

\begin{proof}
{Both claims follow immediately by induction from the construction
of the measure. The projection of the measure
$\nu_{\gamma}$ to the first row (column) is of full support, i.e., 
its support coincides with
$({\bf R}_+)^{\infty}$. In order to check this, it suffices to recall
that the sequence of entries of the first row is a sequence of
i.i.d.~variables on ${\bf R}_+$ with distribution $\gamma$, 
and hence this sequence is stationary.
Assume that we have proved that the projection is of full support
and the sequence is stationary for the first $k$ rows. Since in the 
above inequality determining the distribution of the entry
$r_{k+1,n}$, both left- and right-hand sides depend only on the 
entries of the previous $k$ rows and this dependence
is invariant under the horizontal shift in $n$ for fixed $k$,
it follows that the left- and right-hand sides of the inequality determine
stationary sequences of random variables. The choice of the entry
$r_{k+1,n}$ reduces to independently choosing a number
uniformly distributed
in the interval between the boundaries of the inequality. 
Thus the elements of the sequence
$\{r_{k+1,n}\}_{n=k+1}^{\infty}$ ($k$ is fixed) are functions
(which do not depend on $n$) of several stationary sequences,
namely the previous rows and a sequence of independent variables;
hence the first $k+1$ rows form a stationary random sequence in
$n$, $n=k+1,k+2, \dots$. Obviously, the shift of the first row
is ergodic; the induction step proceeds as follows: it is easy to see that
the shift acting in the space of $(k+1)$-rows is 
the factor shift of the direct product of the shift acting in the space of
$k$-rows and a Bernoulli shift (that determines the choice of the entries
$r_{k+1,n}$, $n=k+1,k+2, \dots$, in the intervals; see above);
hence it is ergodic. The same considerations imply that 
the projection of the measure is of full support.}
\end{proof}

As we have already mentioned, this method of defining measures
on the set of matrices can be called Markovian, because  
we can also regard the conditional measures defined above as
the transition probabilities of a random field 
if we will regard the random
matrix itself as a random field on the lattice
$({\cal Z}_+)^2$ with boundary condition given by the distribution of the first row
and the first column. The procedure described above is simpler than that given in
\cite{V2}, though they are essentially equivalent. The measure
$\gamma$ on the half-line, as well as the conditional measures, 
may vary in wide ranges; however, these variations are inessential;
only explicit conditions on the random distance matrix
are of importance; they correspond
to the intuitive notion of randomness and guarantee the conditions
of the theorem from Section 3, which we will check with the help of the lemma proved above.

An illustrative description of our construction of the measure is as follows:
given an $n$-point metric space, we add the next 
$(n+1)$th point at random, the vector of distances between the new point
and the old ones (i.e., an admissible vector) having a uniformly
positive probability on the set of admissible vectors.
Thus this procedure can be used for simulating a generic metric
space, as in the Monte-Carlo method. We must warn those who dealt with
experiments with random configurations that metric spaces having an isometric
realization in an Euclidean space form an exponentially
small fraction with respect to our measure as
$n \to\infty$.

Let us supplement the construction with several remarks.

 1. If we impose an additional condition on all distances:
$r_{i,j} \leq b$, then the same procedure (in this case we may assume, 
for example, that 
$\gamma$ is the uniform measure on
$[0,b]$) will give us a random metric space of
diameter $b$.

2. An obvious drawback of our construction is that it is not
invariant with respect to the numbering of points; for example, 
the first point is distinguished: the sequence of distances
from the first point to the other ones is a sequence of
independent variables, but the sequence of distances from the second
(or any subsequent) point to the other ones is no longer a sequence
of independent variables. The situation can be remedied by symmetrization, i.e.,
by considering the weak limit of measures
$$ \bar \nu_{\gamma} = \lim (n!)^{-1}\sum_{g \in S_n} g \nu_{\gamma},$$
where the symmetric group acts by simultaneous permutations
of the rows and columns of a distance matrix (and hence on measures
defined on the cone of distance matrices). It is not difficult to prove
that the weak limit $ \bar \nu_{\gamma}$ exists and, obviously, is invariant
under the action of the infinite symmetric group.
However, this gives rise to a new series of questions, which are partially
discussed in Section 5.

3. One can use the same recipe to construct arbitrary measures on matrices:
first one should determine the distribution of the first
row (column) and then successively determine the conditional
measures of entries. This method applies to constructing
not only random distance matrices, but also measures
on other families of matrices.

\section{Universal distance matrices and the Urysohn space. 
Urysohnness criterion}

\subsection{Definitions and the urysohnness criterion}

The following definition plays a key role in our considerations.

\begin{definition}
{{\rm1.} An infinite distance matrix 
$r=\{r_{i,j}\}_{i,j=1}^{\infty}
\in {\cal R}$ is called a {\it universal distance matrix}
if the following condition holds:

for every $\epsilon >0$, every $n \in \bf N$, and every vector
$a=\{a_i\}_{i=1}^n \in A(p_n(r))$, there exists $m \in \bf N$ such that
$\max_{i=1, \dots, n}|r_{i,m}-a_i| < \epsilon$.

In other words, for every $n \in \bf N$, the set of vectors
$\{\{r_{i,j}\}_{i=1}^n\}_{j=n+1}^{\infty}$ is everywhere dense in
the set of admissible vectors $A(p_n(r))$.

{\rm2.} An infinite proper distance matrix
$r=\{r_{i,j}\}_{i,j=1}^{\infty} \in {\cal R}$ is called an almost
universal distance matrix if for every
$n \in \bf N$, the set of all its submatrices of order $n$ of the form
$\{r_{i_k,i_s}\}_{k,s=1}^n$ 
over all $n$-tuples $\{i_k\}_{k=1}^n \subset \bf N$ 
is dense in the cone ${\cal R}_n$.}
\end{definition}

Denote the set of universal distance matrices
(respectively, almost universal distance matrices) by 
  $\cal M$ (respectively, $\cal M'$). We will prove that
  $\cal M$ is not empty, but first we describe some
  properties of universal matrices.

\begin{lemma}
{Every universal distance matrix is almost universal. There exists
an almost universal nonuniversal distance matrix.}
\end{lemma}

\begin{proof}
 {Let $r$ be a universal matrix; let us 
 choose an arbitrary finite distance matrix
$q \in {\cal R}_n$ and prove that, given a positive number
  $\epsilon$, we can find a set
 $\{i_k\}_{k=1}^n \subset \bf N$ such that $\max_{k,s=1,\dots,
 n}|r_{i_k,i_s} - q_{k,s}|<\epsilon$. Since $r^1=\{r_{1,1}=0\}$,
 we have $A(r^1)={\bf R}_+$ (see Section 2.2); since $r$ is universal, 
the sequence $\{r_{1,n}\}_{n=2}^\infty$ is dense in
${\bf R}_+ $; therefore we can choose a positive integer
$i_1$ such that
 $|r_{1,i_1}-q_{1,2}|<\epsilon$ and then, using the density of columns of 
 length two, which follows from the universality condition,
choose a positive integer
$i_2$ such that $|r_{1,i_2}-q_{1,3}|<
 \epsilon$, $|r_{2,i_2}-q_{2,3}|< \epsilon$; iterating this procedure
for successive columns of the matrix
$q \in {\cal R}_n$, we obtain a desired submatrix of $r$.

There are many examples of almost universal matrices that are not
universal; in particular, the distance matrix of every
countable everywhere dense subset of the universal, but not
homogeneous (see below) Polish space (for example, 
the Banach space $C([0,1])$) yields such an example.}
\end{proof}

\begin{corollary}{\rm(}on $\epsilon$-extension of isometries{\rm)}.
 {Let  $r$ be an infinite distance matrix. The matrix
  $r$ is universal if and only if for every positive
 $\epsilon$, every positive integer $n<N$, and every
finite distance matrix $q \in {\cal R}_N$ of order $N$
whose NW-corner of order $n$ coincides with the NW-corner 
 of the same order of the matrix $r$
{\rm(}i.e., $r_{i,j}=q_{i,j}$, $i,j=1, \dots, n${\rm)}, there exist
positive integers $i_1, \dots, i_N$ with  $i_k=k$, $k= 1, \dots, n $,
such that $\max_{k,s=1, \dots, N}|r_{i_k,i_s}-q_{k,s}|<\epsilon$.

In other words, we can extend the set of the first $n$ positive integers 
with an  $(N-n)$-tuple of positive integers 
 $i_{n+1}, \dots, i_N$ in such a way that the resulting distance  matrix
of order $N$ with indices 
$i_1=1$, $i_2=2$, $\dots$, $i_n=n$, $i_{n+1}, \dots, i_N$
will differ from the distance matrix $q$
{\rm(}with respect to an arbitrary given norm{\rm)}
by less than $\epsilon$.}
\end{corollary}

The property stated in the first claim of the corollary is called
the {\it $\epsilon$-extension of isometries}, 
because it means that an isometry between the first $n$ points
of a finite metric space of cardinality $N$ and a countable metric space
with a given distance matrix can be extended to an
$\epsilon$-isometry of the whole finite metric space.

\begin{proof}
{Assume that $r$ is a universal matrix and $N=n+1$; 
then the assertion follows immediately from the definition of
universality; the desired assertion follows by induction on $N$.
And if the condition of the corollary holds for
$N=n+1$, it follows that the matrix is universal.}
\end{proof}

Another useful reformulation of the notion of universality
and almost universality, which follows from the corollary,
involves group-theoretic terms. Let
$q \in {\cal R}_n$; denote by ${\cal R}^n(q)$
the set of all matrices $r \in \cal R$ such that the NW-corner of $r$ 
coincides with $q$. Consider the group 
$S^n_{\infty}$ of simultaneous permutations of rows and columns
of matrices that leave the first $n$ rows and columns of distance
matrices invariant and hence send
${\cal R}^n(q)$ to itself. The following {\it universality criterion}
follows immediately from definitions.

\begin{statement}
{A matrix $r \in \cal R$ is universal if and only if
for every positive integer $n$, its orbit under the action of the group
$S^n_{\infty}$ is everywhere dense in
${\cal R}^n(r^n)$ in the weak topology {\rm(}here
$r^n$ is the NW-corner of order $n$ of the matrix
$r${\rm)}. A matrix $r\in \cal R$
is almost universal if and only if
for every $n$, its orbit under the action of the whole group
$S_{\infty}$ is everywhere dense in the cone $\cal
R$.}
\end{statement}

The existence of universal distance matrices is not quite obvious;
as we will see, it is equivalent to the existence of the universal Urysohn
space. Urysohn's work
\cite{U1,U,Ur} begins with a direct construction of a concrete 
universal rational distance matrix, though he uses neither this notion, nor
its properties. We slightly simplify Urysohn's proof from
\cite{U} and prove not only the existence of universal matrices, but
also a more general fact:

\begin{theorem}
{The set $\cal M$ of universal matrices is not empty and is an everywhere
dense $G_{\delta}$-set in the cone $\cal R$ (in the weak topology).}
\end{theorem}

Before proving the theorem, note that
the random distance matrix constructed in the previous section
is universal with probability one; this fact will be proved
in the next section; thus, in particular, the
existence of universal matrices will be proved. In this section we will
give a deterministic version of the same inductive construction
of a universal matrix.

\begin{proof}
{Let us first prove the existence of universal distance matrices in the cone
${\cal R}$; we use the inductive method and the notions 
defined in the previous sections.

As the first row of the matrix, take an arbitrary sequence
$\{r_{1,n}\}_{n=1}^{\infty}$
everywhere dense in ${\Bbb R}_+$. Since
$r_{1,1}=0$ and $A(0)$ coincides with
${\Bbb R}_+$, the universality condition for the first row is satisfied.
Assume that we have already constructed
the first $k-1$ rows satisfying the 
desired conditions. Let us define the
$k$th row $\{r_{k,n}\}_{n=k+1}^{\infty}$ elementwise, beginning 
(since the matrix is symmetric) with the entry
$r_{k,k+1}$. By the induction hypothesis, the column vectors
$\{r_{i,j}\}_{i=1}^{k-1}$ for all
$j=k+1,k+2, \dots $ lie in the admissible set $A(r^{k-1})$, where
$r^{k-1}$ is the submatrix of order $k-1$ (the NW-corner of order
$k-1$) of the matrix 
$r$, and the sequence of these columns is dense in
the admissible set $A(r^{k-1})$. We must supplement the columns
with numbers $\{r_{k,j}\}_{j=k+1}^{\infty}$ so that every resulting
column will lie in the set $A(r^k)$ and the whole set of them 
will be dense in $A(r^k)$.
The existence of  such numbers is an immediate consequence of 
the following elementary lemma on the ``lifting'' of a countable
everywhere dense set from the base space of a
fiber bundle to the bundle itself.

\begin{lemma}
{Let $K$  be a closed convex set, in the Euclidean space
${\Bbb R}^k$, represented as a fiber bundle over its projection
$K_0$ to the coordinate hyperplane
${\Bbb R}^{k-1}$, the fiber over a point of 
$K_0$ {\rm(}the preimage of a point under the projection{\rm)} being 
a closed interval of the real axis. For every countable everywhere
dense set $\{y_n\}_{n \in \Bbb \N}$ in $K_0$, there exists a sequence of 
real numbers $\{z_n\}_{n \in \Bbb \N}$ such that every vector
$(y_n,z_n)$ belongs to the cone $K$ {\rm(}in other words,
$z_n$ lies in the interval-fiber over the point
$y_n${\rm)} and the set of all such vectors  is everywhere dense in
$K$.}
\end{lemma}

\begin{proof}
{Indeed, since the fiber over the point
$y_n$ is a nonempty interval
$[a_n,b_n]$, we can choose a finite or countable everywhere dense set
$\{z_{i,n}\}_i$ in each of these intervals. 
Obviously, we can define the points
$z_n=z_{i(n),n}$ corresponding to the points $y_n$ 
by induction in such a way that the countable
set $\{(y_n, z_n)\}_i$ will be everywhere dense in
$K$. We omit standard details.}
\end{proof}

We apply this lemma as follows: let $K=A(r^{k})$; then its projection $K_0$ is
$A(r^{k-1})$, and the fiber over a point $\{r_{i,j}\}_{i=1,\dots, k-1} \in
A(r^{k-1})$ is the interval of possible values for the matrix entry
 $r_{k,j}$, $j=k+1,k+2, \dots $. The existence of universal matrices 
 is proved.
 
Since the universality of a distance matrix is a property that is preserved under
simultaneous finite permutations of rows and columns, 
as well as under the NW-shift (see above), and moreover
it is preserved under a change of finitely many matrix entries, it follows that
the set $\cal M$ contains, together with each matrix, all its permutations
and NW-shifts. Since
universality implies almost universality, it follows that the 
$S_{\bf N}$-orbit of
a universal matrix $r$ 
is everywhere dense in
$\cal R$ in the weak topology. Thus
$\cal M$ is everywhere dense in $\cal R$.

Finally, the following formula is an immediate consequence of the definition
of universality and shows that the set
$\cal M$ of universal matrices is an everywhere dense
$G_{\delta}$-set:
$${\cal M} =\bigcap_{k \in {\bf N}}\bigcap_{n \in {\bf N}}
\bigcap_{a \in A(r^n) \bigcap {\bf Q}^{n^2}} \bigcup_{m \in {\bf N},
m>n}\{r\in {\cal R}:\max_{i=1,\dots,
n}|r_{i,m}-a_i|<\frac{1}{k}\}.$$}
\end{proof}

Fix a universal proper distance matrix 
$r$ and equip the set of all positive integers 
$\bf N$ with the metric
$r$. Denote by $({\cal
U}_r, \rho_r)$ the completion of the metric space
$({\bf N},r)$ with respect to the metric
$r$. It is obviously a Polish space.

\begin{lemma}
The distance matrix of every countable everywhere dense subset
$\{u_i\}$ of the space ${\cal U}_r$ is universal.
\end{lemma}

\begin{proof}
{Identify the set of positive integers
${\bf N}$ with the set $\{x_i\} \subset {\cal U}_r$.
By definition, $\rho(x_i,x_j)=r_{i,j}$. Since $r$ is a universal matrix,
it follows that for every $n$, the following property
holds (where $\Cl$ stands for the closure):
$\Cl(\cup_{j>n}\{\{\rho(x_i,x_j)\}_{i=1}^n\}=A(p_n(r))$. Since
the set $\{u_i\}$ is everywhere dense in $({\cal U}_r, \rho_r)$, 
we may replace the last inequality with the following one:
$\Cl(\cup_{j>n}\{\{\rho(x_i,u_j)\}_{i=1}^n\}=A(p_n(r))$. But since
$\{x_i\}$ is also everywhere dense in $(U_r, \rho_r)$, we can write
$\Cl(\cup_{j>n}\{\{\rho(u_i,u_j)\}_{i=1}^n\}=A(p_n(r'))$, where $r'$
is the distance matrix of the  sequence
$\{u_i\}$.}
\end{proof}

Thus the universality of the distance matrix of a countable 
dense subset is a property of the Polish space --- it does or
does not hold for all such subsets simultaneously.
We will see that the space $({\cal U}_r, \rho_r)$ 
is the so-called universal Urysohn space (which is
defined below) and the universality of a distance matrix
is a necessary and sufficient condition that guarantees that 
the completion of the set of positive integers with respect to 
the corresponding metric is isometric to the Urysohn space.

\subsection{Universal Urysohn space and universal matrices}

Now let us define the Urysohn space. In his last paper
\cite{U}, which was published already after his tragic death,
Pavel Urysohn (1898--1924) suggested a construction of the universal
Polish space, which is now called the ``Urysohn space.''
This paper was an answer to the question, posed by M.~Fr\'echet,
on the existence of the universal Banach space.
Later Banach and Mazur gave a direct answer to Fr\'echet's question
by showing that there exists a universal Banach space
(for example, $C([0,1])$); their result is a simple consequence
of a theorem due to Hausdorff. In the book
\cite{Ba}, Urysohn's paper is mentioned, but without 
indicating the most important properties of the Urysohn space:
homogeneity and uniqueness up to isometry. But precisely  
these properties impart a fundamental character to this space. 
We unify and slightly generalize the main theorems of the paper
\cite{U} in the following statement.

\begin{theorem}{\rm(}Urysohn {\rm\cite{U1,U})}

{There exists a Polish {\rm(}= complete separable metric{\rm)}
space $\cal U$ with the following properties:

{\rm1)} {\rm(}Universality.{\rm)} For every Polish space
$X$, there exists an isometric embedding of
$X$ into $\cal U$.

{\rm2)} {\rm(}Homogeneity.{\rm)} For any two isometric finite subsets
$A=(a_1, \dots, a_m)$ and $B=(b_1, \dots, b_m)$ of the space $\cal U$,
there exists an isometry $J$ of the whole space
$\cal U$ that sends $A$ to
$B$: $JA=B$.

{\rm3)} {\rm(}Uniqueness.{\rm)} Any two Polish spaces that 
satisfy conditions {\rm1)} and {\rm2)} are isometric.}
\end{theorem}

Under the conditions of the theorem, property 2) can be strengthened:
finite subsets can be replaced with compact ones; thus the group
of isometries of the space acts transitively on isometric compact subsets; 
however, compact subsets in this
statement cannot be replaced with closed ones.

Condition 2) can be strengthened in a similar way and formulated
as a condition of extension of a given isometry of two compact subsets
to the whole space; more precisely:

2') {\it Let $F$ be a finite subset in $\cal U$,  and let $i$
be an isometric embedding of $F$ into a finite metric space
$F'$. Then there exists an isometric embedding 
$I$ of the space $F'$ into $\cal U$ such that the product $I \cdot i$
is the identity mapping of the subset $F$.}

Below we will prove the homogeneity in this strengthened form.
In this condition we can also replace finite sets
$F$ and $F'$ with compact ones.

Note that there are many distinct (nonisometric) universal, but not
homogeneous Polish spaces (for example, the Banach space
$C([0,1])$). The above corollary shows that the principal difference
between such universal spaces and the Urysohn space is as follows:
any Polish space can be isometrically embedded in any universal space;
but in the case of the Urysohn space, one can do much more:
given a fixed finite (and even compact) subset of the Urysohn space,
the images of the points of the Polish space under the embedding can have 
any prescribed distances to the points of this subset. 
Similarly, the homogeneity property holds, for example,
in Euclidean and Hilbert spaces. But the two properties simultaneously
determine the space uniquely up to isometry.

The main result of this section is the following theorem, which
includes the previous one.

\begin{theorem}
{\rm1.} The completion  $({\cal U}_r, \rho_r)$ of the space of positive
integers $({\bf N},r)$ with respect to the metric determined by a universal
proper distance matrix $r$ satisfies properties
{\rm1)},
{\rm2)}, and {\rm2')} of the previous theorem, i.e., is the Urysohn space.

{\rm2. (}Uniqueness.{\rm)} For any two universal proper distance matrices
$r$ and $r'$, the completions of the spaces
$({\bf N},r)$ and
$({\bf N},r')$ are isometric. Thus the isometric type of the space
$({\cal U}_r, \rho_r)$ does not depend on the choice of the universal 
matrix $r$. The universality is a necessary and
sufficient condition on the distance matrix of any countable
everywhere dense subset of the Urysohn space.
\end{theorem}

The proof of Theorem 3 given below partially reproduces and simplifies
the arguments of Urysohn's paper. Urysohn did not use infinite
distance  matrices. We will essentially follow our paper
\cite{V2}.

\begin{proof}
{Assume that $r=\{r(i,j)\}_{i,j=1}^{\infty} \in {\cal
R}$ is a proper universal distance matrix (it is convenient to write
$r(i,j)$ instead of $r_{i,j}$) and ${\cal U}_r$ is the completion
of the countable metric space $({\bf N},r)$. We denote the corresponding
metric on ${\cal U}_r$ by $\rho_r$,
sometimes omitting the index $r$.

1. First of all, let us prove that the metric space
$(\cal U,\rho)$ is universal in the sense of property
1) from Theorem 2 and that it is homogeneous in the sense of property 2)
from the same theorem.

Let $(Y,q)$ be an arbitrary Polish space. In order to prove that there 
exists an isometric embedding of
$(Y,q)$ into $(\cal
U,\rho)$, it suffices to prove that there exists  an isometric embedding 
of an arbitrary countable dense subset
$\{y_n\}_1^{\infty}$ of the space $(Y,q)$ into $({\cal U},\rho)$. 
Thus we must prove that for every infinite proper distance matrix
$q=\{q(i,j)\} \in \cal R$, there exists a dense subset
$\{u_i\} \subset \cal U$ with distance matrix equal to
$q$. This in turn means that we must construct 
a set of fundamental sequences in the space
$({\bf N}, r)$, say $N_i=\{n^{(m)}_i\}_{m=i}^{\infty}
\subset {\bf N}$, $i=1,2, \dots$, such that
\begin{eqnarray*}
\lim_{m \to {\infty}} r(n^{(m)}_i, n^{(m)}_j)&=&q(i,j),\quad i,j=m,m+1,
\dots, \quad\mbox{and} \\
\lim_{m,k \to \infty}
r(n^{(m)}_i,n^{(k)}_i)&=&0\quad\mbox{for all } i.
\end{eqnarray*}

 The convergence of the sequence
 $N_i=\{n^{(m)}_i\}_{m=i}^{\infty}$ in the space
$({\cal U},\rho)$ as $m \to \infty$ to a point
$u_i\in {\cal U}$, $i=1,2 \dots $, follows from the fundamentality of the
sequence, i.e., from the second equation; the first equation implies that
the distance matrix of the limiting points
$\{u_i\}$ coincides with $q$. 

Let us now construct the sequences
$N_i$, $i=1,2, \dots$. We will construct the desired sequences
$\{N_i\}_{i=1}^{\infty} \subset \bf N$ by induction.

Choose an arbitrary point $n^{(1)}_1 \in \bf N$ and assume that
for a given $m>1$, we have already determined the finite fragments
$L_k=\{n^{(k)}_i\}_{i=1}^k \subset {\bf N}$ of the first  $m$ sets
$\{N_i\}_{i=1}^m$ for $k=1, 2, \dots, m$ such that
$\max_{i,j=1,\dots, k} |r(n^{(k)}_i, n^{(k)}_j)-q(i,j)|=\delta_k <
2^{-k}$, $k=1,2, \dots, m$, and the sets $L_k$ 
are pairwise disjoint.

Our construction of the set
$L_{m+1}$ depends only on the set
$L_m$; thus for simplicity we can renumber the sets
$L_m$: $n^{(m)}_i=i$, $i=1, \dots, m$.

Let us construct a new set $L_{m+1}=\{n^{(m+1)}_i\}_{i=1}^{m+1}
\subset \bf N$ with desired properties as follows.
Consider a finite metric space
$(V,d)$ consisting of $2m+1$
points $y_1, \dots, y_m; z_1, \dots, z_m, z_{m+1}$ with distances $
d(y_i,y_j)= r_{i,j}$, $i,j=1, \dots, m$, $d(z_i,z_j)= q(i,j)$,
$i,j=1,\dots, m+1$; $d(y_i,z_j)=q(i,j)$, $i=1, \dots, m$, $j=1,\dots, m+1$; $i
\ne j$,  $d(y_i, z_i)= \delta_m$, $i=1, \dots, m$, for some
$\delta_m$. It is easy to check that these distances are well defined.
Denote the distance matrix of the space
$(V,d)$ by $q_m$. Apply Corollary 2 (on $\epsilon$-extension of isometries)
and extend the set $L_m=\{1, 2, \dots, m\}$ with a new set
$L_{m+1}$ consisting of  $m+1$ positive integers
$\{n^{(m+1)}_i\}_{i=m+1}^{2m+1} \subset \bf N$ so that the distance matrix
of the new space $L_{m+1}$ will differ from the NW-corner of order  
$m+1$ of the matrix $q$ by at most $\delta_m$, which is less than
$2^{-(m+1)}$: $\max_{i,j} |r(n^{(m+1)}_i, n^{(m+1)}_j)
- q_m(i,j)| =\delta_{m+1}<2^{-(m+1)}$ {\sloppy} (recall that the
NW-corners of order $m$ of the matrices $q_m$ and $r$ coincide
by construction). We see that for every $i$, the sequence
$\{n^{(m_i)}\}_{m=i}^{\infty}$ is fundamental and $\lim_{m \to
\infty} r(n^{(m)}_i,n^{(m)}_j)=q(i,j)$. Thus we have proved that every Polish
space can be isometrically embedded in
$({\cal U}_r,\rho_r)$.

Now we can substantially strengthen Corollary 2:

\begin{corollary} {\rm(}on extension of isometries{\rm)}.
{The space $({\cal U}_r,\rho_r)$ enjoys the following property:
for every finite set $A=\{a_i\}_{i=1}^n \in {\cal U}_r$
and every distance matrix $q$ of order $N$, $N>n$, with NW-corner of order $n$
equal to the distance matrix $\{\rho(a_i,a_j)\}_{i,j=1}^n$,
there exist points $a_{n+1}, \dots, a_N$ such that the distance matrix
of the whole set $\{a_i\}_{i=1}^N$ equals $q$.}
\end{corollary}

The proof parallels that of Corollary 2 and uses
the arguments given above.

Let us continue the proof of Theorem 3.

2. In order to prove the homogeneity in the strong form
2'), take two arbitrary finite $n$-point isometric
subsets $A=\{a_i\}_{i=1}^n$  and $B=\{b_i\}_{i=1}^n$ in  $({\cal
U}_r,\rho_r)$ and choose an isometry between them by fixing the order
of points corresponding to the isometry. Let us construct two
isometric ordered countable subsets
$C$ and $D$, each of them dense in
${\cal U}$, such that the first  $n$ points of $C$ 
coincide with the points of
$A$ in the chosen ordering, and the first points of
$D$ coincide with the points of
$B$.  First we fix a countable everywhere dense
subset $F$ in $({\cal U}_r,\rho_r)$ such that $F \cap A =F \cap
B=\emptyset $ and represent it as an increasing union:
$F=\cup F_n$. Let $C_1=A \cup F_1$ and find a set
$D_1=B \cup F'_1$ such that the isometry of the subsets 
$A$ and $B$ extends to $F_1$ and
$F'_1$. Thus $D_1$ is isometric to $C_1$. 
This can be done by Corollary 4 (on extension of isometries).

Then choose $D_2=D_1 \cup F_2$ and $C_2=C_1 \cup F'_2$ 
and again extend the isometry from the subsets $D_1$ and $C_1$ 
to the whole sets. Thus we have constructed an isometry
between $D_2$ and $C_2$, etc. The alternating process yields two
everywhere dense countable isometric sets
$\cup C_i$ and $\cup D_i$, the isometry between them extending 
the isometry between the subsets
$A$ and $B$. The method we have used is well known and is called 
in the literature
the ``back and forth'' method; it was also used in Urysohn's paper.

 3. Uniqueness. Let
$r$ and $r'$ be two universal proper distance matrices,
and let  $({\cal U}_r, \rho_r)$ and $({\cal
U}_r', \rho_r')$ be the corresponding completions of the set
of positive integers.

Let us construct in these spaces two countable everywhere dense
subsets $F_1$ and $F_2$, respectively, so that the isometry between them 
extends to an isometry of the whole spaces. Denote by    
$\{x_i\}$ and $\{u_i\}$ the everywhere dense subsets of the spaces
$({\cal U}_r, \rho_r)$ and
$({\cal U}_r', \rho_r')$ that generate the matrices 
$r$ and $r'$, respectively. Now we repeat the same argument as in the
proof of the first claim of the theorem. We begin with a finite set of points
$\{x_i\}_{i=1}^{n_1}$ in $({\cal U}_r, \rho_r)$ and supplement it with a set of points
$\{u'_i\}_{i=1}^{m_1} \subset {\cal U}_r$ with the same distance matrix as
that of the set 
$\{u_i\}_{i=1}^{m_1}$; this can be done by the universality of the space
$({\cal U}_r, \rho_r)$ (property 1), which is already proved. Now supplement
the set $\{u_i\}_{i=1}^{m_2}$ ($m_2 > m_1$) in the space
${\cal U}_r'$ with a set of points $\{x'_i\}_{i=1}^{n_2}$
($n_2 > n_1$) so that the distance matrix of the subset  
$\{u_i\}_{i=1}^{m_1} \cup \{x'_i\}_{i=1}^{n_1}$ of the set
$\{u_i\}_{i=1}^{m_2} \cup \{x'_i\}_{i=1}^{n_2}$ will coincide with
the distance matrix of the set $\{u'_i\}_{i=1}^{m_1} \cup
\{x_i\}_{i=1}^{n_1}$, etc; continuing this process to infinity, 
we obtain the desired two sets: 
$\{x_i\}_{i=1}^{\infty} \cup \{u'_i\}_{i=1}^{\infty} \subset {\cal
U}_r$ and $\{u_i\}_{i=1}^{\infty} \cup
\{x'_i\}_{i=1}^{\infty} \subset {\cal U}_r'$, 
which are everywhere dense in their spaces and isometric. Thus
we have completed the proof of the theorem.}
\end{proof}

Theorems 1 and 3 yield the following remarkable fact:

\begin{corollary}
{The generic distance matrix is a universal matrix; and hence the generic
Polish space {\rm(}in the sense of our model, i.e., the cone $\cal R${\rm)} 
is the Urysohn space $\cal U$.}
\end{corollary}

In our terms, P.~Urysohn first chooses a rational
distance matrix, which determines the countable universal (noncomplete) 
metric space in the class of metrics with rational values, and its
completion is the universal space
$\cal U$. Urysohn (\cite{U}) observes that there exists a universal space
of diameter 1 (or any given diameter). If we similarly define
universal distance matrices with elements from the interval
$[0,1]$, then the corresponding completion is precisely the universal space
for Polish spaces of diameter 1, and all theorems of this
section remain valid after trivial modifications. The question,
posed in \cite{U}, whether there exist noncomplete universal (homogeneous)
spaces was answered in the affirmative in
\cite{Ka}: it turns out that one can remove special subsets from
the Urysohn space with preserving homogeneity, universality, and
density of the remaining set. Katetov's construction of the 
universal space is slightly different from the original one and from our
constructions; a close version was suggested by M.~Gromov in
\cite{G}\footnote{It seems that the papers by Katetov
\cite{Ka} and Sierpinski \cite{Ser} were the first ones, after many years,
to deal with the Urysohn space; textbooks on measure-theoretic topology
rarely make mention of it.}.

Recently V.~Uspensky \cite{Us} (see also \cite{Us1}),
using a theorem due to Torunczyk
\cite{To}, proved, answering my question, that the Urysohn space
is homeomorphic to the infinite-dimensional Hilbert space
(and hence to any separable Banach space). Although this result
is hardly useful for constructing a realization of the Urysohn space,
because the homeomorphism between the Urysohn space and
the Hilbert space is highly nonconstructive and, which is most 
important, has no uniformness or smoothness, 
nevertheless it is of theoretical interest, because it shows
that the Urysohn space can be regarded as a standard space
in the infinite-dimensional
topology, perhaps more natural than the Hilbert space. But the problem
of constructing a
more explicit realization of the Urysohn space, which was mentioned 
already by P.~S.~Alexandrov in his comments to the Russian translation
of the main Urysohn's paper \cite{U}, is still actual: in contrast to,
for example, illustrative models of the universal graph, we have no
nonapproximative models of the Urysohn space.

In \cite{Us}, Uspensky proved the universality of the group $\Iso(\cal U)$ of isometries
of the Urysohn space as a topological group 
in the compactly open topology: every Polish group can be 
continuously isomorphically embedded into
$\Iso(\cal U)$; the key point is the proof of the existence of an
isometric embedding of an arbitrary Polish space into
$\cal U$ such that its group of isometries is naturally isomorphically embedded
into the group of isometries of the Urysohn space. V.~Pestov
\cite{Pe} proved that the group
$\Iso(\cal U)$ has a fixed point property, i.e., every continuous action
of this group on an arbitrary compact space has a fixed point.
Among more special and not quite obvious properties of the Urysohn space,
we mention that is is not isometric to its direct powers with natural metrics
(F.~Petrov) and that replacing its metric with any nonlinear concave function
of the original metric also leads to spaces that are not isometric to
the Urysohn space (i.e., are not universal). Thus the main open problems
are as follows: to find a satisfactory model
of the universal space and to describe the structure of the group
of its isometries. In particular, the following special problem 
seems to be rather interesting: whether there exists an everywhere dense
subgroup of the group of isometries that is the inductive limit of
compact groups. A similar combinatorial problem (concerning the group
of automorphisms of the universal graph) is recently solved,
see the Appendix.

\section{The main theorem: universality of almost all distance matrices} 

We will consider arbitrary Borel probability measures on the cone 
$\cal R$, or, in other words, arbitrary random metrics on the
set of positive integers. Since the cone
$\cal R$ is metrizable and separable in the weak topology, i.e., 
becomes a Polish space if we fix some metric compatible with
the weak topology, it follows that the simplex of measures with
the weak topology (topology of convergence on cylinder sets of the cone
$\cal R$) is also a Polish space. This simplex is the inverse limit of 
the simplices of probability measures on the finite-dimensional cones
${\cal R}_n$.

Recall that in Section 2.3 we have constructed a family of measures
$\nu_{\gamma}$ on the cone
${\cal R}_n$, indexed by an arbitrary
measure $\gamma$ of full support on the half-line.

\begin{theorem}
{The measures $\nu_{\gamma}$ constructed in Section {\rm2.3}
are concentrated on the set of universal matrices; in other words, almost
every distance matrix with respect to any of these measures is universal. 
Therefore, the completion of the set of positive integers with respect
to the metric determined by almost every matrix is isometric to the 
Urysohn universal space.}
\end{theorem}

\begin{proof}
{We must check that $\nu_{\gamma}$-almost every
matrix $r=\{r_{i,j}\}$ satisfies the universality condition.
We strengthen this condition as follows:
for every $k$, the empirical distribution of the column
consisting
of the first $k$ entries (for almost every given matrix) is 
of full support, i.e., the support of this empirical distribution
coincides with the set $A(r^k)$. This is indeed a strengthening, since
the universality condition requires only that the column
vectors of length $k$ of a given matrix be dense in
$A(r^k)$ and we will show that they appear with frequency 
corresponding to the given distribution on the set
$A(r^k)$, which is of full support by construction. But, in
view of Lemma 5, the measure  $\nu_\gamma$ is invariant and ergodic with
respect to the shifts of the first  $k$ rows for all
$k$, and thus we can apply the ergodic theorem (= the law of
large numbers), which asserts that the empirical distribution 
coincides with the theoretical distribution for 
almost all matrices, 
and the same lemma implies that the original
distribution is of full support for every $k$.}
\end{proof}

Thus we have proved that $\nu_{\gamma}$-almost every
distance matrix generates, as the completion of the 
set of positive integers with respect
to the corresponding metric, the Urysohn space. These measures should be
thought of as the most ``dispersed'' over the set of distance matrices;
hence this theorem justifies the assertion that the random space
is universal with probability one. Another paradoxical property of these
measures will be considered in Section 5.

However, we can assert much more. Denote by
$\cal V$ the simplex of all Borel probability measures on the cone 
$\cal R$ and equip it with the weak topology, i.e., the topology
of convergence on cylinder sets. Note that the set of nondegenerate
(i.e., positive on nonempty open sets) measures is 
an everywhere dense $G_{\delta}$-set in $\cal V$. 
The corollary of Theorems 1 and 3 immediately implies the following
claim on genericity.

\begin{theorem}
{The subset of measures on $\cal V$ concentrated on universal matrices
is an everywhere dense $G_\delta$-set in $\cal V$. Thus, with respect
to the generic measure $\nu$ on $\cal R$,  $\nu$-almost every distance matrix
is universal and hence determines
a metric $r$ on the set of positive integers $\bf N$ 
such that the completion of $\bf N$
with respect to this metric
is the Urysohn space.}
\end{theorem}

This assertion follows immediately from the general fact that the set of all
measures defined on a Polish space and concentrated on a fixed
everywhere dense  $G_\delta$-subset is itself an everywhere dense
$G_\delta$-subset of the space of all measures in the weak topology.
Thus we have another confirmation of the assertion that 
the Urysohn space is generic, now in the probabilistic sense: {\it the random
countable metric space is isometric to an everywhere dense subset
of the Urysohn space}, or, in other words, the completion of the random
countable metric space is the Urysohn space with probability one;
the ``randomness'' is understood with respect to any measure from a
$G_{\delta}$-subset. However, this assertion, unlike the theorem proved above,
yields no interesting concrete examples of such measures (as the measures
$\nu_{\gamma}$). In the next section we will study measures on the
cone of distance matrices in more detail.

\section{Matrix distributions and a generalization of Kolmogorov's problem of 
extension of measures}

\subsection{Invariants of metric triples, and the uniqueness theorem}

We will consider metric spaces with measure and random metrics on
the set of positive integers. Assume that
$(X,\rho,\mu)$ is a Polish space with metric
$\rho$ and Borel probability measure
$\mu$. We say that $(X,\rho,\mu)$ is a 
{\it metric triple} (in
\cite{G}, the term ``mm-space'' = metric-measure space is used; other
terms are ``probability metric space,'' ``Gromov triple'').
Two metric triples
$(X_1,\rho_1, \mu_1)$ and $(X_2, \rho_2, \mu_2)$
are {\it isomorphic} if there exists a measure-preserving  isometry $V$:
$$\rho_2(Vx,Vy)=\rho_2(x,y),\quad  V \mu_1=\mu_2.$$

As we have already mentioned, the classification of (noncompact) Polish
spaces is a ``wild'' problem. It is surprising that the classification of
metric triples is a ``tame'' problem and has a reasonable answer in terms
of the action of the infinite group of permutations
$S_{\infty}$ and $S_{\infty}$-invariant measures on
$\cal R$.

Given a metric triple  $T=(X,\rho, \mu)$, consider its infinite product
with the Bernoulli measure $(X^{\bf N},\mu^{\bf N})$ and define a mapping
$F: X^{\bf N} \to {\cal R}$ as follows:
$$ F_T(\{x_i \}_{i=1}^\infty) =\{\rho(x_i,x_j)\}_{i,j=1}^\infty
 \in {\cal R}.$$
The $F_T$-image of the measure  $\mu^{\bf N}$, which we denote by $D_T$,
is called the {\it matrix distribution of the metric triple $T$}:
$F_T \mu^{\infty} \equiv D_T$.

The group $S_{\infty}$ of all finite permutations of positive integers
(= the infinite symmetric group $S_{\bf N}$) acts both on
$\bf M_N(R)$ and on the cone $\cal R$ of distance matrices as the group
of simultaneous permutations of rows and columns of matrices.

\begin{lemma}
{The measure $D_T$ is a Borel measure on
$\cal R$, invariant and ergodic with respect to the action of the infinite
symmetric group and the action of the NW-shift 
{\rm(}the simultaneous shift of an infinite matrix in the vertical
(up) and horizontal (to the left) direction:
$(NW(r))_{i,j}=r_{i+1,j+1}; i,j=1,2, \dots${\rm)}.}
\end{lemma}

\begin{proof}
{All claims follow from the analogous properties of the measure
$\mu^{\infty}$, which is invariant under the shift and
permutations of coordinates, and from the fact that the mapping
$F_T$ commutes with the action of the shift and permutations.}
\end{proof}

We say that a measure on a metric space is {\it nondegenerate}
if every nonempty open set has a positive measure.

\begin{theorem}
{Two metric triples $T_1=(X_1,\rho_1, \mu_1)$ and
$T_2=(X_2,\rho_2,\mu_2)$ with nondegenerate measures are equivalent
if and only if their matrix distributions coincide as measures on the cone
  $\cal R$: $D_{T_1}=D_{T_2}$.}
\end{theorem}

\begin{proof}
{It is obvious that the matrix distributions of equivalent
triples coincide: if there is a measure-preserving  isometry
$V:X_1 \to X_2$ between
$T_1$ and $T_2$, then its infinite power
$V^{\infty}$ preserves the Bernoulli measure:
$V^{\infty}(\mu_1^{\infty})= \mu_2^{\infty}$; hence, since
$F_{T_2}X_2^{\infty}=F_{T_2}(V^{\infty}X_1^{\infty}$),
the images of these measures coincide:
$D_{T_2}=D_{T_1}$.

Now assume that
$D_{T_2}=D_{T_1}=D$. Then $D$-almost all distance matrices
$r$ are the images, under the mappings 
$F_{T_1}$ and $F_{T_2}$, of some sequences, say
$r_{i,j}=\rho_1(x_i,x_j)=\rho_2(y_i,y_j)$; but this means that
identifying $x_i \in X_1$ and $y_i\in X_2$ for all $i$ is an isometry
$V$ between these countable sets. Since the measures are nondegenerate,
these sequences are everywhere dense with probability one; hence
the isometry extends to the whole image and preimage.
Let us prove that this isometry preserves the measure. This is the key point:
by the ergodic theorem, 
$\mu_1$-almost all sequences
$\{x_i\}$ and $\mu_2$-almost all sequences
$\{y_i\}$ are uniformly distributed on
$X_1$ and $X_2$, respectively. This means that our sequences 
can be chosen so that the 
$\mu_1$-measure of every ball of rational radius centered at the points
of our sequences 
$B^l(x_i) \equiv \{z\in X_1:\rho_1(x_i,z)<l\}$ will equal 
$$
\mu_1(B^l(x_i)) =\lim_{n \to \infty}\frac{1}{n} \sum_{k=1}^n
1_{[0,l]}(\rho_1(x_k,x_i))
$$
(countably many conditions). But since $V$ is isometric
(because $r_{i,j}=\rho_1(x_i,x_j)=\rho_2(y_i,y_j)$, see above), the same expression is the
$\mu_2$-measure of the ball $B^l(y_i)
\equiv \{u\in X_2:\rho_2(y_i,u)<l\}$:
 $$\mu_2(B^l(y_i))
=\lim_{n \to \infty}\frac{1}{n} \sum_{k=1}^n
1_{[0,l]}(\rho_2(y_k,y_i))= \mu_1(B^l(x_i)). $$

Finally, since the measures are nondegenerate and hence, as mentioned
above, the sequences
$\{x_i\}$ and $\{y_i\}$ are everywhere dense each in its space,
it follows that if the values of two (Borel) measures on a countable set 
of balls with rational radii coincides, then the $V$-image
of the measure in the first space coincides with the measure in
the second space.}
\end{proof}

\begin{corollary}
{The matrix distribution is a complete invariant of equivalence
classes {\rm(}up to measure-preserving isometry{\rm)}
of metric triples with nondegenerate measures.}
\end{corollary}

We call this result the ``uniqueness theorem,'' because it
asserts that a metric triple with a given matrix distribution
is unique up to equivalence. This theorem, in another
form and under the name of ``Reconstruction
Theorem,'' was originally proved in the book \cite[p.~117--123]{G}.
The theorem was stated in terms of finite-dimensional distributions
of the measure that we call the matrix distribution,
and its proof involved analytical techniques. The ``ergodic'' proof
given above was the answer to the question posed by Gromov in
1997: how one can avoid complicated arguments. This proof was also
given in \cite{V1} and cited in the book
\cite[p.~122--123]{G}. The same book invites the reader to compare
both proofs and explain how the ergodic theorem replaces analytic arguments
(such as the Weierstrass theorem and the  method of moments). 
The explanation is as follows:
the ergodic theorem allows one to use, instead of
approximative techniques, an infinite (limiting) object
(infinite orbits, measures on the limiting space, etc.) and its
properties (for example, equidistribution) that
cannot even be formulated for finite objects. In our case, considering
infinite matrices and the cone
$\cal R$ with an invariant measure allows us to reduce the problem
to studying an ergodic action of the infinite symmetric group.
In \cite{V2,V3} we use this ergodic machinery to obtain a much more
general result --- the classification of measurable function of several 
arguments; above we have considered a special case of this problem:
a metric is a function of two arguments on a metric space with measure.

\medskip

\subsection{Existence theorem}

By definition, the matrix distribution of a metric triple
$T=(X,\rho, \mu)$ with a nondegenerate measure
$\mu$ is a measure $D_T$ on the cone $\cal R$.  
In this section we consider random matrices (or measures on
$\cal R$) that can be matrix distributions; in other words,
we want to describe those
distributions on distance matrices that arise as the matrix distributions
of a sequence of independently chosen points
$\{x_i\}$ of some metric space
$(X,\rho)$ distributed according to some measure
$\mu$ on this space. To obtain such a
characterization is also
necessary in order to assert that the classification problem is
indeed ``smooth,'' i.e., that the set of invariants is described explicitly.
We will show that the set of matrix distributions is  a Borel subset
of the space of all probability measures on the cone
$\cal R$.

As mentioned above (Lemma 7), every measure $D_T$ must be
invariant and ergodic under the action of the symmetric group
and NW-shift. But this condition is not sufficient, and there are
additional conditions. We give necessary and sufficient conditions
on the measure 
(see also \cite{V3}); below we present counterexamples that show that
these conditions are substantial.

\begin{theorem}{\rm(}Existence of a metric triple with a given matrix
distribution{\rm)}.

{Let $D$ be a Borel probability measure on the cone of distance matrices
$\cal R$ that is invariant and ergodic under the action of the 
infinite symmetric group.

{\rm1)} The following condition is necessary and sufficient for the measure $D$
to be the matrix distribution of some metric triple
$T=(X,\rho, \mu)$:

for every $\epsilon>0$, there exists an integer 
$N=N(\epsilon)$ such that
\begin{equation}
D\{r=\{r_{i,j}\} \in {\cal R}: \lim_{n \to \infty}
\frac{|\{j:1\leq\ j \leq n, \min_{1 \leq i \leq N}
 r_{i,j} <\epsilon \}|}{n}> 1- \epsilon \}>1-\epsilon.
\end{equation}

{\rm2)} The following stronger condition
is necessary and sufficient for the measure $D$
to be the matrix distribution of some metric triple
$T=(X,\rho, \mu)$ with a compact metric space
$(X, \rho)$:

for every  $\epsilon >0$, there exists an integer 
$N=N(\epsilon)$ such that
\begin{equation}
D\{r=\{r_{i,j}\} \in {\cal R}:  \min_{1\leq
i \leq N} r_{i,j} <\epsilon\mbox{ for all } j>N\}>1-\epsilon.
\end{equation}}
\end{theorem}

Here we will only sketch the proof of this theorem (see
\cite{V2} for details).

 A. Necessity. In the case of a compact space the necessity is obvious:
condition (5) means that a sufficiently long sequence of independent
points of the space distributed according to the measure
$\mu$ contains an
$\epsilon$-net of the compact metric space for every
$\epsilon$. In the general case, the
necessity of condition (4) 
follows automatically from the following
well-known property of Borel measures in
Polish spaces: there is a sigma-compact set of full measure 
(the so-called ``regularity'' of measures); therefore for every
$\epsilon > 0$, there exists a compact set of measure
$>1-\epsilon$. Hence, by the countable additivity of the measure,
for every $\epsilon >0$, there exists finitely many points such that
the measure of the union of the 
$\epsilon$-balls centered at these points is greater than
$1-\epsilon$; using the ergodic theorem, we can assert that
the condition inside the braces in (4) is satisfied for
distance matrices from a set of measure greater than
$1-\epsilon$.

B. Sufficiency. Now assume that $D$ is an invariant and ergodic
measure on the cone $\cal R$ satisfying condition (4). Choose an
arbitrary $D$-generic distance matrix --- one can check that 
$D$-almost all matrices satisfy the desired conditions --- and
construct a complete metric space by completing the set
of positive integers with respect to the corresponding metric; the set
of positive integers forms an everywhere dense family in this space.
Using the classical law of large numbers, one can easily
find the values of the measure on the set algebra generated by the balls
of arbitrary radii centered at the chosen points, and condition (4)
allows one to prove that this measure extends to a countably
additive Borel measure. Then one checks that the matrix distribution
of the obtained metric triple coincides with the measure $D$. 

\smallskip\noindent{\bf Remark}
{The structure of the necessary and sufficient conditions on the measure
in the above theorem shows that matrix distributions form a Borel
subset in the set of all Borel probability measures on the cone 
$\cal R$. Thus the classification of metric triples is a smooth
(Borel) classification, because the space of types is a Borel set of
matrix distributions in the set of Borel probability measures
on the cone of distance matrices. The same conclusion holds for a more
general problem of classification of measurable functions of several 
arguments, see \cite{V3}.}
\smallskip

Condition (4) can be replaced with another condition from the paper
\cite{V3}, the so-called  ``simplicity'' of an 
$S_{\infty}$-invariant (or symmetric) measure, 
which is stated below. This condition guarantees that a 
measure $D$ is the matrix distribution of a measurable function
of two variables, and in view of the uniqueness theorem, 
this suffices
for our purposes. It remains to check that the matrix distribution of the
constructed metric triple coincides with the original measure.

\subsection{Uniform distribution of sequences in metric spaces and
a generalization of Kolmogorov's problem on extension of measures}

Let us turn to a more detailed analysis of measures on the cone
of distance matrices and properties of individual matrices.

Let $(X,\rho)$ be a metric space, and let
$x=\{x_i\}_{i=1}^{\infty}$ be an everywhere dense sequence of points in $X$.
We say that this sequence is
{\it regular} if for every finite $n$ and every open subset
$C \subset R^n$, the limit
 $$\lim_{n \to \infty}\frac{1}{n} \sum_{k=1}^n
1_{C}(\{\rho(x_i,y_k)\}_{i=1}^n)$$
exists. In particular, for every positive integer $i$ and every positive
$r$, the limits
$$\mu(D_r(x_i))\equiv \lim_{n \to \infty}\frac{1}{n} \sum_{k=1}^n
1_{[0,r)}(\rho(x_i,y_k)\})$$
exist.

The latter condition is a specialization of the former one;
it means, in particular, that one can find the ``empirical''
measures of open balls 
$D_r(x_i)$ of radius $r$ centered at the points
$x_i$ of our sequence; and the more general first condition allows one
to find the joint distribution of the distances to an arbitrary
finite number of points of our sequence, i.e., to find the measures
of intersections and unions of balls, etc.; in other words, this
condition allows one to {\it determine a finitely additive measure
$\mu_x$ on the set algebra of the space $X$ generated by the open balls centered
at the points of the sequence}. One may say that a regular sequence
is uniformly distributed with respect to the measure
$\mu_x$ determined by it. If we were given a Borel probability
(countably additive) measure $\bar \mu$ on the space
$(X,\rho)$, then the classical law of large numbers would imply that 
a realization of the sequence of i.i.d.~points of $X$
distributed according to the measure
$\mu$ is regular with probability one, and the above conditions show
how one can recover the measure. But if the measure is not given {\it a priori},
the question arises: whether for every regular sequence
$x=\{x_i\}$,
$i=1,2,\dots$, the finitely additive measure  
$\mu_x$ extends to a Borel probability measure
$\bar \mu$ on the space
$(X,\rho)$? In other words, whether every regular sequence is
uniformly distributed with respect to some
probability measure?

This question is similar to the questions that arise
in connection with the classical Kolmogorov theorem on extension
of measures and its generalizations: in that setting, given 
a measure defined
on the algebra of cylinder sets of a vector space, one checks
that it is countably additive on the algebra of cylinder sets,
and then
this measure extends to a true probability measure on the $\sigma$-algebra.
By the Kolmogorov theorem, this extension always exists in the linear space
$R^{\infty}$, but in other spaces this is true not for all
cylinder measures. More precisely, our analogy with the Kolmogorov
theorem is as follows.

Consider the mapping $\Pi \equiv \Pi_x : X \rightarrow {\Bbb
R}_+^{\infty}$ that associates with a point $y$ of the space $X$ the
sequence of its distances
$\{\rho(x_i,y)\}$, $i=1,2, \dots$, to the points of the dense set
$x=\{x_i\}$.
Obviously, this mapping sends 
closed balls centered at $x_i$
to cylinder sets, and the set algebra 
generated by these balls to a subalgebra of the algebra of cylinder sets
of the space ${\Bbb R}_+^{\infty}$. Since the sequence
$x=\{x_i\}_i$ is everywhere dense, the mapping
$\Pi_x$ is monomorphic. 
As mentioned above, the regularity of $x$ allows us
to find the joint distributions of coordinates in the image,
and thus a cylinder measure in
${\Bbb R}_+^{\infty}$; by the Kolmogorov theorem, this measure extends to
a probability (countably additive)
measure $\tau \equiv \tau_x$ in ${\Bbb R}_+^{\infty}$. However, this does not
mean that the mapping
$\Pi_x$ is an isomorphism of measure spaces; it is only a homomorphism
of a space with finitely additive measure to a space with probability measure.
Moreover, a measure on the algebra generated by the balls centered
at the points of an everywhere dense sequence does not 
always extend to a countably additive measure on the $\sigma$-algebra;
namely, the following assertion holds.

\begin{lemma}
{A finitely additive measure
$\mu_x$ extends to a Borel probability measure
$\bar \mu_x$ in the space  $(X,\rho)$ if and only if the 
$\tau_x$-measure of $\Pi_x (X)$ is equal to one.}
\end{lemma}

Condition (4) (or the simplicity condition, see below) guarantees that the
extension exists; in general, if this condition is not satisfied, the measure 
is not countably additive: in the general case, the measure
$\tau_x(\Pi_x (X))$ may be equal to zero, as we will see later.  
In this case it is natural to ask whether we can extend the space 
$(X,\rho)$ so that the measure will equal to one;
I have not a complete answer to this question, but in some cases
such a compactification of the space is possible. However, in this
problem we deal with a ``nonlinear'' version of Kolmogorov's 
extension problem, and there hardly exists an exact analog of the unified space
${\Bbb R}^{\infty}$ in which the Kolmogorov theorem holds
for any cylinder measure.

\smallskip\noindent{\bf Example 1}
A trivial example when there is no countable additivity is as follows. Let
$(X,\rho)$ be a countable set of points with unit pairwise distances,
and let $x$ be simply the sequence of all points of the space $X$. This sequence
is obviously regular, because its distance matrix is the matrix with
units outside the diagonal; hence the measure 
$\tau$ --- the joint distribution of distances --- is the $\delta$-measure
at the point $(1,1,  \dots)\in {\Bbb
R}_+^{\infty}$. And on the space $(X,\rho)$ there is only a 
finitely additive measure $\mu$, which vanishes on finite sets
and equals one on cofinite sets.
\smallskip

Before considering more serious examples, let us formulate the question
in terms of measures on the cone of distance matrices and their
symmetrizations. The definition of regularity immediately implies
the following result.

\begin{lemma}
{It is necessary and sufficient for a matrix
$r \in \cal R$ to be the distance matrix of a regular sequence of some
metric space that there exist the empirical distribution of the rows
{\rm(}columns{\rm)} of $r$ as a probability measure in the space
$\Bbb R^{\infty}$. This measure coincides with the measure
$\tau$ defined above.}
\end{lemma}

We will call such distance matrices {\it regular} and
denote by $\Psi$ the mapping from the set of regular matrices to the set
of probability measures in the space
$\Bbb R^{\infty}$ that associates with a regular matrix the empirical
distribution of its rows.

Without dwelling on details, we note that every regular distance matrix
generates, by the symmetrization 
$$ \bar \mu_{B} = \lim (n!)^{-1}\sum_{g \in S_n} 1_B(gr),$$
an $S_{\infty}$-invariant measure on the cone of distance matrices.
Here $B$ is  a cylinder set in the space of matrices, and  
an element $g$ acts on a matrix $r$ by a simultaneous permutation
of rows and columns.

\begin{theorem}
{The following assertions are equivalent:

{\rm1.} A regular everywhere dense sequence $x=\{x_i\}$
of a Polish space $(X, \rho)$ is uniformly distributed with respect
to some Borel probability measure
$\mu$ {\rm(}of full support{\rm)} on the space $(X, \rho)$.

{\rm2.} The symmetrization of the regular distance matrix $\{\rho (x_i,x_j)\}$
is a matrix distribution.}
\end{theorem}

The matrix distribution from assertion 2 is precisely 
the invariant of the metric triple
$(X, \rho, \mu)$ in the sense of Theorem 6.

Thus we have two equivalent questions:

1. When a regular sequence of points of a Polish space is uniformly distributed
with respect to some Borel (countably additive) probability
measure in this space?

2. When a symmetric ergodic measure on the cone of distance matrices
is the matrix distribution of a metric triple?

The answer to the first question is given by the lemma, and the answer
to the second question is given by condition (4) or
a more general condition of {\it simplicity} of a measure
on the space of matrices:

A symmetric (i.e., $S_\infty$-invariant) measure $\lambda$ 
on the cone of distance matrices
$\cal R$ is called {\it simple} if the mapping
$\Psi$ is an isomorphism of measure spaces 
$(\cal R, \lambda)$ and
 $(\Bbb R^{\infty}, \Psi(\lambda))$; in other words, $\lambda$-almost
every distance matrix is uniquely determined by the empirical distribution
of its rows (columns). The notion of simplicity is introduced in
\cite{V3} and makes sense for arbitrary measures on spaces of matrices.

\begin{theorem}{\rm(}see {\rm\cite{V2}).}
A  symmetric ergodic measure on the cone
$\cal R$ is a matrix distribution if and only if it is simple.
\end{theorem}

Sometimes it is rather difficult to check the simplicity condition;
however, it is obvious, for example, that a product measure on the space
of matrices is not simple.

\smallskip\noindent
{\bf Example 2}
First note that {\it any} symmetric matrix with zeros
on the main diagonal and entries
$r_{i,j}$, $i \ne j$, taking values in the interval
$[1/2,1]$ is a proper distance matrix, because in this case the
triangle inequality holds automatically for any triple of indices. Hence
every probability measure $m$ supported by the interval
$[1/2,1]$ generates the product measure
$m^{\infty}$ concentrated on the cone
$\cal R$; this means that all above-diagonal
entries of the matrix 
are independent with respect to the measure
$m^{\infty}$ and identically distributed with distribution
$m$. Obviously, the measure $m^{\infty}$ is invariant and ergodic under
the action of the group of permutations; if $m$ is not a $\delta$-measure,
then $m^{\infty}$ is a continuous (nondiscrete) measure. The law of
large numbers implies that almost every matrix is regular.
Since nondiagonal elements (pairwise distances) in this example are
bounded away from zero, almost every matrix determines a discrete
metric space, which cannot have a continuous measure. Hence the measure
$m^{\infty}$ is not the matrix distribution of any metric triple.
It is equally obvious that the measure
$m^{\infty}$ is not simple. One easily computes that in this example
the measure $\tau_x(\Pi_x(X))$ vanishes.
\smallskip

However, it is not difficult to prove that this is the {\it unique}
(up to replacing the interval  $[1/2,1]$ with any interval of the form
$[a/2,a]$, $a>0$) class of examples when the entries of the random
distance matrix can be independent random variables. Example 1 given above
is a degenerate case: $m=\delta_1$.

Let us say that a regular sequence $x=\{x_i\}$ in a metric space
is {\it completely regular} if any of the following equivalent
conditions holds:

a) the measure
$\mu_x$ extends to a countably additive measure, and hence the sequence
$x$ is uniformly distributed with respect to this measure;

b) the symmetrization of the distance matrix of the sequence $x$
generates a matrix distribution;

c) the $\tau_x$-measure of the image $\Pi_x(X)$ is equal to one.

A regular sequence such that there is no countably additive measure
with respect to which it is uniformly distributed will be called 
{\it anarchical}. It is easy to prove that in a locally
compact metric space (which means in this context that every closed
ball is compact) there are no anarchical sequences.

The above examples show that in discrete metric spaces, anarchical
sequences are simply sequences consisting of all points of the space.

A much more surprising example is given by the Urysohn space.

\begin{statement}
{\rm1.} Consider the measure $\nu_\gamma$ constructed in Section 
{\rm3}; its symmetrization is not a matrix distribution.

{\rm2.} Almost every matrix with respect to this measure determines
an anarchical sequence in the Urysohn space. Thus this sequence is not
uniformly distributed with respect to any probability
measure in this space.
\end{statement}

This result was discovered after 
O.~Bohigos, Ch.~Schmidt, and E.~Bogomolny at my request had made
experiments on computing the several first distributions of the 
symmetrization of the measures 
$\nu_\gamma$ and their spectra. Their computations showed that the spectrum
of the random distance matrix in this case is governed by the 
Wigner semicircle law and that, for example, the support
of the distribution of the symmetrized matrix entry is bounded away from zero.
This means precisely that the symmetrized measure is not the matrix
distribution of any metric triple with a nondiscrete space --- the result that 
we managed to prove. The underlying reason for this phenomenon is that 
the symmetrization of a measure concentrated on universal matrices
may no longer be concentrated on them. This is the case with
the measure $\nu_\gamma$. Thus the ``most random'' distance matrices
determine no measure on the Urysohn space; they are ``too random''
and in a certain sense correspond to the ``white noise'' in
this space. It is not known 
whether this means that there exists a virtual metric space
(say, some compactification) in which these finitely additive measures
extend to a probability measure. However, this effect
emphasizes the similarity between the Urysohn space and the universal
graph (see the Appendix). In order to construct probability measures in
this space, we need universal matrices of more special form. We will
return to this question elsewhere.

 \subsection{Measure-theoretic space of metric triples}

We sought for invariants of metric triples as invariants of a
measurable function of two variables, i.e., invariants of the metric regarded as
a symmetric function on the direct square of the metric space with 
the product measure. Since a Polish space with a Borel
continuous probability measure is isomorphic, as a measure space, 
to the interval
$[0,1]$ with the Lebesgue measure, the problem reduces to classification
of symmetric measurable functions of two variables. In fact, instead of
the ordinary point of view, when one considers the set of all Borel
measures on a given topological or metric space, we, on the contrary,
{\it fix the structure of a measure space (Lebesgue space)
and consider all measurable (semi)metrics on this space}
(see \cite[\S 6]{V6}). Let us consider this point
of view in more detail.

Assume that $(X,\mu)$ is a Lebesgue space with a finite or
$\sigma$-finite measure $\mu$ (for example, the interval
$[0,1]$ or the line with the Lebesgue measure) and
$S_{\mu}(X)$ is the space of all $\bmod 0$ classes of measurable functions
on $(X,\mu)$; define the cone 
${\cal R}^c \subset S_{\mu}(X)$ of measurable metrics, 
i.e., the cone of classes of 
$\bmod 0$ coinciding symmetric functions of two arguments
$\rho :(X\times X, \mu \times \mu) \to
{\bf R}_+$ satisfying the inequality
$$\rho(x,y)+\rho(y,z) \geq \rho(x,z) \mbox{ for }
(\mu \times \mu \times \mu)\mbox{-almost all } (x,y,z)\in (X
\times X \times X).$$ 
It is also natural to assume that
$(\mu\times \mu)\{(x,y): \rho(x,y)=0\}=0$. In order to single out true
metrics, impose on $\rho$ the following purity condition
(\cite{V3}): the partition of
$(X,\mu)$ into the classes $x \sim y \Leftrightarrow
\rho(x,\cdot)=\rho(y,\cdot) \bmod 0$ is the partition into points
$\bmod 0$.

Functions satisfying all these conditions will be called
{\it almost metrics}. It is essential that 
$\rho$ is not an individual function, but a class of
$\bmod 0$ coinciding functions, hence {\it a priori} it is not obvious
that there exists an individual function (lifting) that determines
the structure of a metric space on $X$ in the literal sense.
If $X$ is a countable or finite set with counting measure, then we
obtain the cone $\cal R$ or ${\cal R}_n$ of distance matrices from 
Section 2. Thus the cone
${\cal R}^c$ is a continuous generalization of the cone
$\cal R$ to the case of a continuous measure.

Now assume that the measure $\mu$ is finite and
$\rho \in {\cal R}^c$ is an almost metric; let a measure
$D_{\rho}$ on the space
$\bf M_{\infty}(\bf R)$ be the matrix distribution of the measurable
function $\rho$ (see the definition in the previous section). It follows
from the ergodic theorem that
$D_{\rho}(r \in {\cal
R}:r_{i,k}+r_{j,k}\geq r_{i,k})=1$  for any $i,j,k \in \bf N$,
and hence $D_{\rho}({\cal R})=1$. Using the characterization
of matrix distributions given in Section 4, we obtain
the following assertion.

\begin{statement}
{The measure $D_{\rho}$ is concentrated on the cone
$\cal R$ {\rm(}i.e., $D_{\rho}({\cal
R})=1${\rm)} and is an ergodic  $S_{\infty}$-invariant measure.
Therefore, every almost metric
$\rho \in {\cal R}^c$ on $(X
\times X, \mu \times \mu)$ determines a $\bmod 0$ metric on the space
$(X,\mu)$.}
\end{statement}

\begin{corollary}
{The set of classes of $\bmod 0$ coinciding measurable semimetrics
on a Lebesgue space with continuous measure coincides with the set
of classes of $\bmod 0$ coinciding metric triples with continuous
{\rm(}finite{\rm)} measure.}
\end{corollary}

\section*{Appendix: examples of other categories where randomness
and universality coincide}

\begin{flushright}
``Consider some total and \\ hence unique copy of
$A$.''

P.~O.~di Bartini \\  Soviet Math. Dokl., vol.~163, no.~4, p.~861--864 (1965).

\end{flushright}
\bigskip

\section*{General setting. Universal and random graphs}

Consider a category in which there is a universally attracting object, i.e.,
an object such that every other object is ``isomorphic'' to its subobject.
Assume that at the same time we are given a definition of the random
object in the category; for example, this is possible if the objects of
the category can be determined by finite or countable collections
of numbers (``structural constants'') satisfying certain conditions,
and random objects are determined by random, in the ordinary sense,
collections of numbers, i.e., by collections of random structural constants
subject to certain conditions. In this case we can pose the following question:
what objects will appear with probability one, and whether they 
will be almost always isomorphic to the universal object.
Of course, the definition of randomness must satisfy certain
requirements, which can be specified following the above
pattern.\footnote{Here it seems especially useful to use Kolmogorov's
complexity-based definition of randomness.}

Within this extremely general setting, which should be
specified (this will be done elsewhere), one may say that the 
case of the category of Polish spaces gives
the positive answer to the question posed above. But it will undoubtedly
remind one of close situations in other categories, and first of all
an analogy with the Erd\"os--R\'enyi theorem
\cite{ER} on random and universal graphs, which asserts that
the random graph is universal with probability one;
see definitions below and in
\cite{Ra,Ca}. This simple theorem is also the simplest case of the above scheme,
because every graph generates a metric space with distance metric 
on the set of vertices, and in the case of the universal graph, the distances
between distinct points assume only two values:
$1$ and $2$. Hence (see the example in Section 2.3) 
the triangle inequality is automatically satisfied.
It is natural to pursue the analogy with metric spaces by using 
distance matrices; but since a graph is usually determined by its adjacency
matrix (which can be uniquely recovered from the distance matrix of the graph
and uniquely determines it), we will also follow this tradition.
More precisely, the random graph in the sense of
\cite{ER} is given by a random adjacency matrix (symmetric
$(0,1)$-matrix) with i.i.d.~entries with probabilities
$(p,1-p)$, where $0<p<1$. Then the distance matrix, which is a
$(1,2)$-matrix, also has a similar distribution. On the other hand, 
the universality of a graph
$\Gamma$ in the sense of \cite{Ra} means that

1) any finite graph is isomorphic to a subgraph of
$\Gamma$;

2) for every $n$, the group of automorphisms of the graph
$\Gamma$ acts transitively on the set of all subgraphs of 
$\Gamma$ of order $n$ isomorphic to an arbitrary given graph.

Using the ``back and forth'' method, one can easily prove that
such a graph is unique up to isomorphism; see, e.g.,
\cite{Ca}. The similarity with the situation considered in this paper
becomes especially strong if we compare the 
universality condition for graphs with the universality criterion 
for metric spaces (Theorem 1). Namely, the following result holds.

Let us say that an infinite symmetric
$(0,1)$-matrix $\{a_{i,j}\}$ with zero diagonal is universal in the class
of $(0,1)$-matrices if for every
$n$ and every word $w_1w_2 \dots w_n$ of length $n$ in the alphabet
$\{0,1\}$, there exists a number 
$k$, $k>n$, such that $a_{i,k}=w_i$,
$i=1,2, \dots, k$, i.e., this word occurs as the beginning of the $k$th column
of the matrix. This definition is the discrete version of the definition
of universality of distance matrices.

\begin{lemma}
{A countable graph is universal if and only if its adjacency matrix is 
universal.}
\end{lemma}

The proof follows immediately from the definition. 
(Cf.~\cite{Ca}, where the universality criterion for graphs is
stated in a slightly different form.) The set of universal matrices is
an everywhere dense $G_{\delta}$-set in the set of all symmetric
$(0,1)$-matrices with zero diagonal in the natural totally disconnected
topology.

But, on the other hand, almost all random graphs in the Erd\"os--R\'enyi
sense obviously satisfy this universality condition, because
all their above-diagonal entries are independent;
moreover, it is easy to deduce from this criterion that there are
many other measures on the space of adjacency matrices with respect
to which almost all graphs are universal; it is sufficient that almost
all matrices contain all binary words as beginnings of their columns.
Such measures also form an everywhere dense
$G_{\delta}$-set in the space of adjacency matrices (or distance
matrices) in the weak topology; Bernoulli measures are only an interval
in this large set.

Thus in the category of graphs there are exact analogs of all
facts proved in this paper for the category of Polish
spaces; one may say that the theory of universal and random
graphs is an exact elementary model of the theory of universal
and random metric spaces. Note that in the literature on graphs, 
the term ``random graph'' traditionally means
``universal graph'';
although this usage of terms
is motivated by the Erd\"os--R\'enyi theorem, it contradicts
the meaning of the notion of randomness, because randomness
can be understood in different ways, but universality is uniquely understood 
and may take place for almost all matrices.

In this connection, let us mention another feature of the theory of
universal graphs, which distinguish them from the Urysohn space.
Since the universal graph is unique up to isomorphism, any two
universal adjacency matrices lie on the same orbit of the infinite
symmetric group, which acts naturally
on the vertices of the graph (or by simultaneous
permutations of the rows and columns of matrices). Hence the group of all
permutations acts transitively on the set of universal matrices, and
we obtain another example of the paradoxical situation first observed
by Kolmogorov (see \cite{V7}): a group that acts transitively has 
a continuum of pairwise singular invariant ergodic measures.

The group of all automorphisms of the universal graph, as well as the
group of all isometries of the Urysohn space, is of extreme interest;
we will return to this subject later. Now we only mention that the recent
paper \cite{Mc} gave the positive answer to 
the question on existence of a dense locally finite
subgroup in the group of automorphisms of the universal
graph. (Compare with the similar question
for the Urysohn space at the end of Section 3.)

The considered example is a very special case
(``toy model'') of the universality of metric spaces studied in this paper;
in this case we have only one additional condition:
the metric assumes only two nonzero values:
$1$ and $2$. Indeed, consider the set of proper distance matrices
whose nondiagonal entries equal
$1$ or $2$. In this case the triangle inequality is automatically satisfied
for all symmetric ($1,2$)-matrices. One can define the notion
of universality of such matrices by analogy with
the definition from Section 3:
every finite sequence of 1's and 2's occurs as the beginning of some
column of the distance matrix. On the other hand, every graph determines a
(geodesic) metric on the set of vertices.

 \begin{lemma}
The set of universal distance $(1,2)$-matrices coincides with the 
set of distance matrices of universal graphs.
 \end{lemma}

Thus the universal graph, regarded as a metric space, is
the universal (and homogeneous) space in the class of
countable metric spaces with metrics assuming two nonzero values:
$1$ and $2$.

In a similar way one can consider universal metric spaces with
metrics assuming finitely many values.

\section*{Several further examples and comments}

Let us briefly mention examples of other categories in which
a similar situation takes place.

\medskip

1. Universal and random simplices.

\medskip

Consider the category of separable compact simplices. It is easy to
construct a simplex whose set of extreme points is everywhere dense ---
the Poulsen simplex {\cal P}. It is less known that this simplex
is universal: every compact simplex can be affinely continuously
embedded in {\cal P} as a face, and any two isomorphic faces can
be sent to each other by an affine isomorphism of the whole simplex
(\cite{Lu,LO}).  Note that the set of extreme points of the universal simplex
in the induced topology is homeomorphic to the infinite-dimensional
Hilbert space (\cite{Lu,LO}).

Let us now define the notion of the random simplex by analogy with
the random metric space. For this purpose, define an equipped simplex as
a simplex with a distinguished countable (or finite) set of
extreme points dense in the whole set of extreme points. An equipped simplex
has a natural continuous embedding into the space
${\Bbb R}^{\infty}$, and thus it is determined by the matrix
whose rows are the vectors corresponding to the images of the distinguished
set of extreme points. The simplex can be obtained as the closed convex
hull in ${\Bbb R}^{\infty}$ of the set of these vectors. It is not
difficult to describe conditions on matrices that arise in this way;
they will again form a closed convex cone in the space of
infinite matrices. Now it is clear how to define the random simplex:
it is determined by a random matrix, or a measure in the space
of matrices concentrated on the above cone.
As in the case of metric spaces, one can give examples of
distinguished Bernoulli-type measures.

A series of assertions similar to the main theorems of this paper is
as follows:

1) The Poulsen universal simplices form an everywhere dense $G_{\delta}$-set
 in the set of all equipped simplices.

2) With respect to the distinguished measures, almost all
simplices are universal.

3) The set of measures
enjoying property 2) is an everywhere dense 
$G_{\delta}$-set in  the set of all measures on the cone
with the weak topology of measures.

Thus we can justify the thesis that a random simplex is universal.

\medskip
 2. Random Banach spaces.

\medskip

Following the same pattern, one can construct the universal 
$L^1$-type space (Gurarij space), and equipped spaces of this kind
again can be determined
by matrices --- see \cite{LO}. Random matrices will correspond to
random $L^1$-type Banach spaces, and with respect to appropriate
measures, almost every space will be universal. 

\medskip
3. Consider metric spaces $(X,r)$ in which the triangle inequality
holds in the following strengthened form:
$$r(x,y)^p + r(y,z)^p \geq r(x,z)^p, \quad p \in [1,\infty).
$$
For $p=1$, this is the ordinary triangle inequality. In the limiting case 
$p=\infty$, the triangle inequality turns into the ultrametric inequality 
$$ \max (r(x,y),r(y,z)) \geq r(x,z).$$
For brevity, spaces that satisfy the inequality with parameter
$p \in [1,\infty]$ will be called $p$-metric.

\begin{statement}
{{\rm1)} For all $p \in [1,\infty)$, there exists the universal
$p$-metric space in the class of Polish spaces. The random
{\rm(}in the above sense{\rm)} $p$-metric space is universal 
with probability one.

{\rm2)} There is no universal ultrametric
{\rm(}$p=\infty${\rm)} space.}
\end{statement}

The proof of this proposition exactly reproduces the constructions of
our paper: in the same way one defines the cone of matrices, universal matrices, 
measures, etc. It is interesting that in the case $p=\infty$ 
the sets of admissible vectors
(see Section 2) degenerate into an amoeba, which
implies assertion 2). The nonexistence of the universal ultrametric space
was also proved in \cite{Bo}. If we fix a countable set of real
numbers (for example, the set of rational numbers with denominator equal
to a power of some fixed simple number) and consider
ultrametric spaces such that the distances between their points assume
values only from this set, then in the class of such spaces the universal
space may exist (this is the case in the example cited
in the parentheses above).

At the same time, random ultrametric spaces can be defined in the same
way, but most probable they will not be isometric to each other.
This question will be considered elsewhere. Besides, here we do not
consider examples of categories of algebraic origin:
universal and random groups, algebras, lattices, etc.

\medskip
4. Let us consider ordinary metric spaces but require that all distances
take values in some subset
$L\in \Bbb R$ of real numbers. For example, if this subset is 
$\{0,1,2\}$, then, as can be easily seen, in the set of such metric
spaces there exists a universal (with respect to finite spaces)
homogeneous space: this is exactly the Rad\'o graph
$\Gamma$ discussed above, regarded as a metric space (with geodesic
metric); recall that any two
vertices of the universal graph
either are connected by an edge, or have a vertex that is
connected with both of them.

Let us consider countable metric spaces with {\it rational distances}.
Here there also exists a universal metric space; 
Urysohn's work begins exactly with the construction of this space, and its completion 
is the universal metric space
$\cal U$. However, it seems to be an important object by itself.
In this case
the notion of a universal matrix is defined as above, 
the urysohnness criterion (see Section 3) has the same form, 
and one can prove the existence of the rational
universal matrix in exactly the same way as above. 
The most important object is the group of isometries
of this space. It is a closed subgroup of the group
$S^{\infty}$ of all permutations of the countable set (in the weak topology),
and the problem of existence of a locally finite everywhere dense subgroup of
this group is of principal interest. See above for a similar question on the 
group of isometries of the Urysohn space.

\smallskip
Translated by N.~V.~Tsilevich.

\bigskip

\hbox{A.~Vershik} \hbox{St.~Petersburg Department of Steklov
Institute of Mathematics}
\hbox{vershik@pdmi.ras.ru}

\end{document}